\documentclass{amsart}
\usepackage{amsmath,amsfonts,amssymb,amsthm}
\usepackage{amsfonts}
\usepackage[mathscr]{eucal}
\usepackage{mathrsfs}
\usepackage{framed}
\usepackage[colorlinks=false,linktocpage=true]{hyperref}
\usepackage{tocvsec2}
\usepackage{latexsym}
\usepackage{color}
\newtheorem{thm}{Theorem}[section]
\newtheorem{cor}{Corollary}[section]
\newtheorem{lem}{Lemma}[section]
\newtheorem{prop}{Proposition}[section]

\theoremstyle{definition}
\newtheorem{defn}{Definition}[section]

\newtheorem{rem}{Remark}[section]
\newtheorem{notn}{Notation}[section]
\newcommand{\thmref}[1]{Theorem~\ref{#1}}

\newcommand{\lemref}[1]{Lemma~\ref{#1}}

\newcommand{\notnref}[1]{Notation~\ref{#1}}
\def\qed{\quad\vcenter{\hrule\hbox{\vrule height.6em\kern.6em\vrule}\hrule}}
\newenvironment{pf}{{\bigskip\textit{\newline Proof.}\quad}}{$\qed$\bigskip\newline}
\newenvironment{pf*}[1]{{\bigskip\textit{\newline#1.}\quad}}{$\qed$\bigskip\newline}
\numberwithin{equation}{section}
\font\tenscrpt=eusm10
\font\sevenscrpt=eusm10 scaled 700
\font\fivescrpt=eusm10 scaled 500
\newfam\eusmfam
\textfont\eusmfam=\tenscrpt
\scriptfont\eusmfam=\sevenscrpt
\scriptscriptfont\eusmfam=\fivescrpt

\def\s{\bold s}

\def\t{\bold t}

\def\Wpsx{W^{x}(\s)}
\def\Wx{W^{x}}

\def\p{\partial}

\def\eqdef{\overset{\triangle}{=}}

\def\R{\Bbb{R}}

\def\Nn{\Bbb{N}_{n}}

\def\Rpn{{\Bbb{R}}_+^n}

\def\bRpn{\partial{\Bbb{R}}_+^n}
\def\intRpn{\overset{\hspace{1mm}\circ\quad}{\Rpn}}

\def\Rpnmo{{\Bbb{R}}_+^{n-1}}
\def\Rp{{\Bbb{R}}_+}
\def\Rd{{\Bbb{R}}^{d}}
\def\P{\Bbb{P}}

\def\EP{{\Bbb E}_{\Bbb P}}

\def\sU{\mathscr U}

\def\sF{{\mathscr F}}

\def\OFP{(\Omega,\sF,\P)}

\def\tnj{{{\bold t}_{\not j}}}

\def\snj{{{\bold s}_{\not j}}}

\def\utnj{u_{{\bold t}_{\not j}}}

\def\Wsnjrx{W^{x}_{\snj}(r_{j})}
\def\Wsnjzx{W^{x}_{\snj}(0)}

\def\WxBn{{\mathbb W}^{x}_{B^{(1)},\ldots,B^{(n)}}}
\def\WxBnt{{\mathbb W}^{x}_{B^{(1)},\ldots,B^{(n)}}(\t)}
\def\WxEnt{{\mathbb W}^{x}_{\Lambda^{(1)},\ldots,\Lambda^{(n)}}(\t)}

\def\ds{\displaystyle}

\def\eqdef{:=}

\def\sm{\setminus}

\def\lap{\Delta}

\def\df#1#2{\ds{\frac{#1}{#2}}}

\def\lbl#1{\label{#1}}

\def\pa{\partial}

\def\lqv{\left<}
\def\rqv{\right>}

\def\lab{\left|}
\def\rab{\right|}
\def\lpa{\left(}
\def\rpa{\right)}
\def\lbk{\left[}
\def\rbk{\right]}
\def\lbr{\left\{}
\def\rbr{\right\}}
\def\bdf{\begin{defn}}
\def\edf{\end{defn}}
\def\bcr{\begin{cor}}
\def\ecr{\end{cor}}
\def\bnt{\begin{notn}}
\def\ent{\end{notn}}
\def\brm{\begin{rem}}
\def\erm{\end{rem}}
\def\blm{\begin{lem}}
\def\elm{\end{lem}}
\def\bpf{\begin{pf}}
\def\bpfs{\begin{pf*}}
\def\epf{\end{pf}}
\def\epfs{\end{pf*}}

\def\beq{\begin{equation}}
\def\beqs{\begin{equation*}}
\def\eeq{\end{equation}}
\def\eeqs{\end{equation*}}
\def\bsp{\begin{split}}
\def\esp{\end{split}}
\def\bc{\begin{cases}}
\def\ec{\end{cases}}
\def\bt{\begin{tabular}}
\def\et{\end{tabular}}

\def\bthm{\begin{thm}}
\def\ethm{\end{thm}}
\def\bsthm{\begin{sthm}}
\def\esthm{\end{sthm}}
\def\bpr{\begin{prop}}
\def\epr{\end{prop}}
\def\bfr{\begin{framed}}
\def\efr{\end{framed}}
\def\bsh{\begin{shaded}}
\def\esh{\end{shaded}}
\def\bcm{\iffalse}
\def\babs{\begin{abstract}}
\def\eabs{\end{abstract}}
\def\ben{\begin{enumerate}}
\def\rencomalp{\renewcommand{\labelenumi}{(\alph{enumi})}}
\def\rencomrom{\renewcommand{\labelenumi}{(\roman{enumi})}}
\def\een{\end{enumerate}}

\def\KBstxy{K^{\textsc{BS}(n,d)}_{\t;x,y}}
\def\KBssxy{K^{\textsc{BS}(n,d)}_{\s;x,y}}
\def\KBsnjrjxy{K^{\textsc{BS}(n,d)}_{\snj,r_{j};x,y}}
\def\KBsnjrhojxy{K^{\textsc{BS}(n,d)}_{\snj,\rho_{j};x,y}}
\def\KBmitzs{K^{\textsc{BM}}_{t_{i};0,s_{i}}}
\def\KBmjtzs{K^{\textsc{BM}}_{t_{j};0,s_{j}}}

\def\KEmitzs{K^{{\Lambda,\beta}}_{t_{i};0,s_{i}}}
\def\KEmjtzs{K^{{\Lambda},\beta}_{t_{j};0,s_{j}}}

\def\KEmtzx{K^{{\Lambda},\beta}_{t;0,x}}
\def\KEmhalftzx{K^{{\Lambda},\frac12}_{t;0,x}}
\def\KEmhalftzxsqrt2{K^{{\Lambda},\frac12}_{t;0,\sqrt2 x}}
\def\KabBmtzxsqrt2{K^{{|B|}}_{t;0,\frac{x}{\sqrt{2}}}}
\def\KabBmtzx{K^{{|B|}}_{t;0,x}}

\def\P{{\mathbb P}}

\def\EP{{\mathbb E}_{\P}}

\def\R{{\mathbb R}}
\def\Rd{{\mathbb R}^d}

\def\S{{\mathbb S}}
\def\Rp{{\mathbb R}_+}

\def\R{{\mathbb R}}

\def\sUj{\sU^{j}}
\def\sVj{\mathscr{V}^{(j)}}
\newcommand{\ceil}[1]{\lceil{#1}\rceil}
\def\pd{\partial}
\begin{document}
\title[BTBS, ISLTBS, and interacting fractional and $\Delta^{\nu}$ PDE{\scriptsize s} systems]{Interacting time-fractional and $\Delta^{\nu}$ PDE{\scriptsize s} systems via Brownian-time and Inverse-stable-L\'evy-time Brownian sheets}
\author{Hassan Allouba}
\address{Department of Mathematical Sciences, Kent State University, Kent,
Ohio 44242}
\email{allouba@math.kent.edu}
\author{Erkan Nane}
\address{Department of Mathematics and Statistics, Auburn University, Alabama 36849}
\email{ezn0001@auburn.edu}
\date{\today}
\begin{abstract}  Lately, many phenomena in both applied and abstract mathematics and related disciplines have been expressed in terms of high order and fractional PDEs.  Recently, Allouba introduced the Brownian-time Brownian sheet (BTBS) and connected it to a new system of fourth order interacting PDEs.  The interaction in this multiparameter BTBS-PDEs connection is novel, leads to an intimately-connected linear system variant of the celebrated Kuramoto-Sivashinsky PDE, and is not shared with its one-time-parameter counterpart.  It also means that these PDEs systems are to be solved for a family of functions, a feature exhibited in well known fluids dynamics models.  On the other hand, the memory-preserving interaction between the PDE solution and the initial data is common to both the single and the multi parameter Brownian-time PDEs.   Here, we introduce a new---even in the one parameter case---proof that combines stochastic analysis with analysis and fractional calculus to simultaneously link BTBS to a new system of temporally half-derivative interacting PDEs as well as to the fourth order system proved earlier and differently by Allouba.  We then introduce a general class of random fields we call inverse-stable-L\'evy-time Brownian sheets (ISLTBSs), and we link them to $\beta$-fractional-time-derivative systems of interacting PDEs for $0<\beta<1$.  When $\beta=1/\nu$, $\nu\in\lbr2,3,\ldots\rbr$, our proof also connects an ISLTBS to a system of memory-preserving  $\nu$-Laplacian interacting PDEs.  Memory is expressed via a sum of temporally-scaled $k$-Laplacians of the initial data, $k=1,\ldots, \nu-1$.  Using a Fourier-Laplace-transform-fractional-calculus approach, we give a conditional equivalence result that gives a necessary and sufficient condition for the equivalence between the fractional and the high order systems.  In the one parameter case this condition  automatically holds.
\end{abstract}
\maketitle
\newpage
\tableofcontents
\section{Introduction and statement of results}\setcounter{thm}{-1}
Many phenomena in mathematical physics, fluids dynamics and turbulence models, mathematical finance, and the modern theory of stochastic processes were recently described through fractional and higher order evolution equations (e.g., see \cite{Abtbs}--\cite{AX}, \cite{BanDeB}, \cite{BMN}--\cite{DeB}, and \cite{Fun}--\cite{Podlubny}).   Another feature important in fluids dynamics and other applied fields is the interaction between different quantities of interest (velocity and pressure in the Navier-Stokes example).  In his article \cite{Abtbs}, Allouba introduced $n$-parameter $d$-dimensional Brownian-time Brownian sheet (BTBS)---a Brownian sheet where each ``time'' parameter is replaced with the modulus of an independent  Brownian motion---and connected it to a new system of fourth order interacting PDEs.   The interaction in these PDEs added novel and intricate new features not present in their one parameter ($n=1$) $d$-dimensional ($d\ge1$) Brownian-time-Brownian-motion PDE counterpart---first given in \cite{Abtp1,Abtp2}  (DeBlassie also gave a different proof in \cite{DeB}).  This interaction means, among other things,  that the PDEs are to be solved for a family of $n+1$ functions and not a single function as in the one-time-parameter case.  As in the $n=1$ case, the BTBS-PDEs are memory preserving.  I.e., the solution family for the BTBS fourth order PDEs also interacts with the initial data via a time-scaled Laplacian of the initial function.

In this article, we first introduce a new---even in the one parameter case---proof that judiciously combines It\^o's rule with properties of the Brownian sheet as well as analysis and fractional calculus  to simultaneously link BTBS to a new system of temporally half-derivative interacting PDEs as well as to the fourth order system proved in \cite{Abtbs} by Allouba via a different approach.  Like the fourth order PDEs connection, the interacting feature of this half-derivative PDEs system belongs solely to the multiparameter setting and does not exist in its one-parameter counterpart.  Here, we note that the half-derivative connection to the one-parameter version of BTBS---the Brownian-time Brownian motion---was first noted and established by Allouba and Zheng in \cite{Abtp1} through their half-derivative generator, an implicit equivalence of the Brownian-time Brownian motion fourth order PDE to a half-derivative-in-time fractional PDE.   The fractional PDE connection and its equivalence to the fourth order PDEs was then given explicitly and generalized by Nane,  Meerchaert, Baeumer, Vellaisamy, Orsingher, and  Beghin \cite{Nanesd,MNV09,BMN,BOap09,BOspa09}.
In the second part of the paper, we introduce the class of inverse-stable-L\'evy-time Brownian sheets (ISLTBS's) and, using an adaptation of our proof for the BTBS, we link any given ISLTBS to a $\beta$-fractional-time-derivative system of interacting PDEs for $0<\beta<1$.  When, $\beta=1/\nu$, $\nu\in\lbr2,3,\ldots\rbr$, our proof also connect our ISLTBS to a system of high order interacting PDEs involving the spatial $\nu$-Laplacian term\footnote{Frequently used notations and acronyms can be found in Appendix \ref{B}} $\Delta_{x}^{\nu}$.  Like the BTBS case, the interacting nature of these PDEs is unique to the multiparameter setup and it vanishes in the one-time parameter framework.   The memory-preserving feature for ISLTBS is apparent in its high ($2\nu$) order PDEs via a sum of temporally-scaled spatial $k$-Laplacians of the initial data, $k=1,\ldots, \nu-1$. In all cases, using a Fourier-Laplace-transform-fractional-calculus approach, we give a conditional equivalence result that gives a necessary and sufficient condition for the equivalence between the fractional and the high order systems.   In the one parameter case this condition is automatically satisfied.  In a separate article \cite{AN}, we treat the case of stable-L\'evy-time fractional Brownian sheet connection to PDEs systems.

We now introduce the setup and state our main results.   It is both instructive and motivating to first introduce and treat the Brownian-time Brownian sheet case.
\subsection{Recalling the Brownian-time Brownian sheet and its fourth order interacting PDEs system}
Let $B^{(1)},\ldots,B^{(n)}$ be $n$ independent copies of a standard one-dimensional Brownian motion  starting at $0$ and independent of an $n$-parameter ($n\ge1$) $\Rd$-valued Brownian sheet
$$W^{0}=\lbr W^{0}(\t)=\lpa W^{0}_{1}(\t),\ldots,W^{0}_{d}(\t)\rpa; \t=(t_{1},\ldots,t_{n})\in\Rpn\rbr,$$
``starting'' at $0\in\Rd$ under $\P$---$\P\{W^{0}(\t)=0\}=1$ for $\t\in\bRpn$ (see \notnref{papernot})---all defined on a probability space $\OFP$.  Of course, the Brownian sheet coordinates $\lbr W^{0}_{1}(\t); t\ge0\rbr,\ldots,\lbr W^{0}_{d}(\t); t\ge0\rbr$ are assumed independent. For any $x=\lpa x_{1},\ldots,x_{d}\rpa\in\Rd$, let
\beq\lbl{BSx}
W^{x}:=W^{0}+x=\lbr W^{x}(\t)=\lpa W^{x_{1}}_{1}(\t),\ldots,W^{x_{d}}_{d}(\t)\rpa; \t=(t_{1},\ldots,t_{n})\in\Rpn\rbr.
\eeq
The Brownian sheet $W^{x}$ transition density is given by
\begin{equation}\lbl{bsdensity}
K^{BS(n,d)}_{\s;x,y}=\frac{\exp\bigg(\frac{-|x-y|^2}{2\prod_{i=1}^n s_i}\bigg)}{\bigg(2\pi \prod_{i=1}^n s_i\bigg)^{d/2}}; s_{i}>0,\ i=1,\ldots,n\ \mbox{ and }x,y\in\Rd.
\end{equation}

We define the $n$-parameter $\Rd$-valued Brownian-time Brownian sheet (BTBS) based on $\Wx$ and $B^{(1)},\ldots, B^{(n)}$, and starting at $x\in\Rd$, by
\beq\lbl{BTBSdef}
\bsp
\WxBnt\eqdef \Wx\lpa \lab B^{(1)}(t_{1})\rab,\ldots,\lab B^{(n)}(t_{n})\rab\rpa;\mbox{ } \t\in\Rpn.
\end{split}
\eeq
Clearly, $\P\lbk\WxBnt=x\rbk=1$ for $\t\in\bRpn$.
Our first main result gives a nonlinear fourth order interacting PDE and a corresponding linear system of fourth order interacting PDEs that are solved by running the BTBS in \eqref{BTBSdef}.  Before stating our main result, it is helpful to adopt some simplifying notational conventions.
\begin{notn}\lbl{papernot} We always use the notation $\Nn:=\lbr1,\ldots,n\rbr$; and we denote by $\intRpn$ and $\bRpn$ the interior and boundary of $\Rpn$, respectively.  We will alternate freely between $u_{\t}$, $u(\t)$, and $u(t_{1},\ldots,t_{n})$ for typsetting convenience and ease of exposition.  Moreover, $\tnj=(t_{i})_{i\in\Nn\setminus\{j\}}\in\Rpnmo$ will denote the $(n-1)$-tuple point in $\Rpnmo$ that is obtained from $\t=(t_{1},\ldots,t_{n})\in\Rpn$ by removing the $j$-th variable, $t_{j}$. The notations $u_{\t}$, $u(\t)$, and $\utnj(t_{j})$ will all mean the $n$-parameter function $u:\Rpn\to\Rd$ evaluated at $\t=(t_{1},\ldots,t_{n})$; and we use $\utnj(t_{j})$ whenever we need to focus on the changes in $u$ as a function of the $j$-th variable $t_{j}$ while holding the rest of the parameters in $\tnj$ fixed (e.g., when we apply It\^o's rule in $t_{j}$).
The same comment applies for the notations $u(\t,x)$ and $u_{\tnj}(t_{j},x)$.  In denoting partial derivatives of any positive integer order $\kappa$ in the variable $t_{i}$, we use $\partial^{\kappa}_{t_{i}}$ and the mixed derivative in $t_{1},\ldots,t_{\kappa}$ by $\partial^{\kappa}_{t_{1},\ldots,t_{\kappa}}$.
\end{notn}
We recall first Allouba's fourth order interacting BTBS-PDEs connection.  Here, we only state the linear system connection (see \cite{Abtbs} for the correspondig nonlinear fourth order interacting PDE).
\bfr
\bthm[Allouba \cite{Abtbs}: BTBS and its fourth order interacting PDEs connections]\lbl{A-BTBSPDE}
Let $\{\WxBnt; \t\in\Rpn\}$ be an $n$-parameter $\Rd$-valued BTBS based on a Brownian sheet $\Wx$ and Brownian motions $\lbr B^{(i)}\rbr_{i=1}^{n}$ and starting on $\bRpn$ at $x\in\Rd$ on $\OFP$.
Let $f:\Rd\to\R$ be bounded and measurable such that all second order partial derivatives $\p^{2}_{x_{k},x_{l}}f$ are bounded and H\"older continuous, with exponent $0<\alpha\le1$, for $1\le k,l\le d$.   If
\beq\label{BTBSexp0}
\bsp
u(\t,x)&=\EP\lbk f\lpa\WxBnt\rpa\rbk,
\\ \mathscr{U}^{(j)}(\t,x)&=\EP\lbk \lpa\prod_{i\in\Nn\setminus\lbr j\rbr}\lab B^{(i)}(t_{i})\rab\rpa^{2}f\lpa\WxBnt\rpa\rbk;\ j\in\Nn,
\end{split}
\eeq
for $(\t,x)\in\Rpn\times\Rd$, then the family $\lbr u,  \sU^{(j)}\rbr_{j\in\Nn}$ is a solution to 
the system of interacting fourth order linear PDEs:
\begin{equation}\label{btbsystem}
\begin{cases}
(a)\ \displaystyle{\p_{t_{j}} u}=\sqrt{\frac{\prod_{i\in\Nn\setminus\lbr j\rbr}t_i}{2^{4-n}t_j\pi^{n}}}\lap_{x} f+\frac18\lap_{x}^{2}\sU^{(j)};& \t\in\intRpn,x\in\Rd,\\
(b)\ u(\t,x)=f(x);&\t\in\bRpn,x\in\Rd,\\
(c)\ \sU^{(j)}(\t,x)=0;&\tnj\in\p\Rpnmo, x\in\Rd,\\
(d)\ \sU^{(j)}(\t,x)=\lbk \prod_{i\in\Nn\setminus\lbr j\rbr}t_{i}\rbk f(x);&\ t_{j}=0, x\in\Rd,
 \end{cases}
\end{equation}
for $ j\in\Nn.$
\ethm
\efr
The intriguing interaction (or coupling) in \thmref{A-BTBSPDE} between $u$ and $\sU^{(j)}$, $j=1,\ldots,n$, shows that the PDEs in \eqref{btbsystem} are nontrivial intricate generalizations of the Brownian-time Brownian motion PDE connection in \cite{Abtp1,Abtp2} (the case $n=1$), given by
\begin{equation}\label{BTBMPDE}
\begin{cases}
 \partial_{t}u(t,x)=\dfrac{\Delta_{x} f(x)}{\sqrt{8\pi t}}+\dfrac18\Delta_{x}^2u(t,x);& t>0,\,x\in\Rd,
\cr u(0,x)=f(x); & x\in\Rd,
\end{cases}
\end{equation}
in which $u=\sU^{(j)}$.  This coupling phenomenon in  the $n>1$ is caused by the interaction between the Brownian-times $\lab B^{(1)}\rab,\ldots,\lab B^{(n)}\rab$ and the variance of the outer Brownian sheet $\Wx$ through $\pa_{t_{j}}\mbox{Var}\lpa\Wx(t)\rpa=\prod_{i\in\Nn\sm\{j\}}t_{i}$.   As with its one-parameter Brownian-time process (BTP) counterpart, the BTBS is memory preserving as is indicated by the inclusion of the Laplacian of the initial function $\lap_{x}f(x)$ in the BTBS PDEs in \eqref{btbsystem}.  This means that the PDE's dynamics in \eqref{btbsystem} depend on $f(x)$ for $t_{j}>0$, and the role of $f$ carries over beyond its more customary starting time---$t_{j}=0$---in more traditional PDEs.    This preservation of the effect of the initial data $f(x)  $ is a manifestation of the non-Markovian nature of the BTBS in the time parameter $t_{j}$. The effect of the initial data, through $\Delta_{x} f$, tapers off as $t_{j}\to\infty$ at the rate of $\sqrt{({\prod_{i\in\Nn\setminus\lbr j\rbr}t_i}){\lpa 2^{4-n}t_j\pi^{n}\rpa^{-1}}}$.  Also, just as BTPs are not classical (not semimartingales, not Markovian, and not Gaussian), BTBS is not a classical random field.

\subsection{The new BTBS time-fractional interacting PDEs system}
The Caputo fractional derivative of order $0<\beta<1$, defined by
\begin{equation}\label{CaputoDef}
\partial_t^\beta u(t,x)=\frac{1}{\Gamma(1-\beta)}\int_0^t \partial_s u(s,x)\frac{ds}{(t-s)^\beta},
\end{equation}
whenever $u$ is $\mathrm{C}^{1}$ in the time parameter.  One of the main advantages of the Caputo fractional derivative is that initial conditions are expressed in terms of initial values of integer order derivatives. Thus, the Caputo fractional derivative is well suited for the Laplace transform techniques and to properly handle initial values \cite{Caputo}.   It has been widely used to solve ordinary differential equations that involve a fractional time derivative  \cite{GorenfloSurvey,Podlubny}.

Our first main result for the BTBS provides its new link to a temporal half-derivative PDEs system.
\bfr
\bthm[BTBS and its interacting fractional PDEs system]\lbl{AN-BTBSPDE}
Assume that $\{\WxBnt; \t\in\Rpn\}$ is an $n$-parameter $\Rd$-valued BTBS based on a Brownian sheet $\Wx$ and Brownian motions $\lbr B^{(i)}\rbr_{i=1}^{n}$ and starting on $\bRpn$ at $x\in\Rd$ on $\OFP$.
Assume the same conditions of \thmref{A-BTBSPDE} on $f:\Rd\to\R$.   If $u$ is as in \thmref{A-BTBSPDE} and
\beq\label{ISSTBSexp0-half}
\mathscr{V}^{(j)}(\t,x)=\EP\lbk \lpa\prod_{i\in\Nn\setminus\lbr j\rbr} |B^{(i)}(t_{i})|\rpa f\lpa\WxBnt\rpa\rbk;\ j\in\Nn,
\eeq
for $(\t,x)\in\Rpn\times\Rd$, then the family $\lbr u,  \mathscr{V}^{(j)}\rbr_{j\in\Nn}$ is a solution to the system
 \begin{equation}\label{isstbsfrcsystem-half}
\begin{cases}
(a)\ \displaystyle\p_{t_j}^{\frac12} u=\frac{1}{\sqrt{8}}\lap_{x}\mathscr{V}^{(j)};& \t\in\intRpn,x\in\Rd,\\
(b)\ u(\t,x)=f(x);&\t\in\bRpn,x\in\Rd,\\
(c)\ \mathscr{V}^{(j)}(\t,x)=0;&\tnj\in\p\Rpnmo, x\in\Rd,\\
(d)\ \mathscr{V}^{(j)}(\t,x)=\lbk \prod_{i\in\Nn\setminus\lbr j\rbr}\sqrt{\frac{2}{\pi}}t_{i}^{\frac12}\rbk f(x);&\ t_{j}=0, x\in\Rd,
 \end{cases}
\end{equation}
for $ j\in\Nn$.  In particular, if $n=1$, $u=\sVj$ and \eqref{ISSTBSexp0-half} reduces to
 \begin{equation}\label{BTBM-half}
\begin{cases}
(a)\ \displaystyle\p_{t}^{\frac12} u=\frac{1}{\sqrt{8}}\lap_{x}u;& t\in(0,\infty),x\in\Rd,\\
(b)\ u(0,x)=f(x);&x\in\Rd,
 \end{cases}
\end{equation}
\ethm
\efr
Our proof of this result is stochastic analytic  in flavor, and it is illuminating since it unifies the fourth order and fractional PDEs in one new equation (see Subsection \ref{newsapf} below).  In doing so, it
\ben\rencomrom
\item simultaneously establishes \thmref{A-BTBSPDE} and \thmref{AN-BTBSPDE} for the BTBS, and in the process it gives another proof to \thmref{A-BTBSPDE} different from that given in \cite{Abtbs} (this also
provides a novel unifying way to concurrently prove the fourth order and fractional PDE connections to the one parameter ($n=1$) Brownian-time Brownian motion case); and it
\item leads to a revealing proof (new even in the one parameter $n=1$ case) in the more general case of inverse L\'evy-time Brownian sheet case that links the ISLTBS to two new interacting PDEs systems: one fractional in time and the other of order $2\nu$, $\nu\in\lbr2,3,4,\ldots\rbr$.
\een

It is now important to highlight the significant difference between the fourth-order-fractional PDEs connection in the one-parameter Brownian-time Brownian motion case and that in the multiparameter Brownian-time Brownian sheet case here.    Unlike \eqref{BTBMPDE} and \eqref{BTBM-half}, which are both single PDEs in one function $u$, the $n>1$ cases in \eqref{btbsystem} and \eqref{isstbsfrcsystem-half} are systems of interacting PDEs involving two \emph{different} families of functions $\lbr u,\sUj\rbr_{j\in\Nn}$ and $\lbr u,\sVj\rbr_{j\in\Nn}$.  This means that an equivalence between the fourth order system \eqref{btbsystem} and the half-derivative system \eqref{isstbsfrcsystem-half} may only be established by relating $\sUj$ and $\sVj$ for $j\in\Nn$.  This is the content of our second main result given in the following subsection.

The Caputo half-derivative feature of absorbing the initial condition $u(0,x)$ in PDEs as compared to the first derivative (compare \eqref{btbsystem} with \eqref{isstbsfrcsystem-half}) can be seen in several different ways.  The first way is through our stochastic analytic proof which simultaneously proves \thmref{A-BTBSPDE} and \thmref{AN-BTBSPDE}.  The second way is through its Laplace transform $s^\beta \tilde u(s,x)-s^{\beta-1} u(0,x)$ which incorporates the initial value as well.
\subsection{Conditional equivalence of the BTBS fourth order and time-fractional interacting  PDEs systems}
Our second main theorem for the BTBS furnishes the condition needed for equivalence to hold between the fourth order and the fractional systems in \thmref{A-BTBSPDE} and \thmref{AN-BTBSPDE}.
\bfr

 \begin{thm}[Conditional equivalence of the fractional and the fourth-order interacting PDEs systems]\lbl{equiv}
 Assume that $f$ satisfies the conditions of \thmref{A-BTBSPDE} and that the pair $\lbr\sUj,\sVj\rbr$ satisfies
\begin{equation}\label{equivalence-cond}
\begin{cases}
(a)\ \displaystyle\sqrt{8}\lap_x(\partial_{t_j}^{1/2}\mathscr{V}^{(j)})=\lap_x^2\mathscr{U}^{(j)}; & \t\in\intRpn,x\in\Rd,\\
(c)\ \sU^{(j)}(\t,x)=\lbk \prod_{i\in\Nn\setminus\lbr j\rbr}t_{i}\rbk f(x);&\ t_{j}=0, x\in\Rd,\\
(e)\ \mathscr{V}^{(j)}(\t,x)=\lbk \prod_{i\in\Nn\setminus\lbr j\rbr}\sqrt{\frac{2}{\pi}}t_{i}^{\frac12}\rbk f(x);&\ t_{j}=0, x\in\Rd,\\
 (d)\ \mathscr{V}^{(j)}(\t,x)=\sU^{(j)}(\t,x)=0;&\tnj\in\p\Rpnmo, x\in\Rd,
 \end{cases}
 \end{equation}
 for some $j^{*}\in\Nn$ (all the derivatives exist and are continuous).  Then \eqref{isstbsfrcsystem-half} is satisfied for $j^{*}$ by $\lbr u, \mathscr{V}^{(j^{*})}\rbr$  iff \eqref{btbsystem} is satisfied for $ j^{*}$  by  $\lbr u,  \sU^{(j^{*})}\rbr$.   Conversely, if $\eqref{btbsystem}$ and $\eqref{isstbsfrcsystem-half}$ are satisfied by $\lbr u,  \sU_{\nu}^{(j^{*})}\rbr$ and $\lbr u, \mathscr{V}^{(j^{*})}\rbr$ for $j^{*}\in\Nn$, then \eqref{equivalence-cond} holds for $j^{*}$.  In particular, the fourth order system $\lbr\eqref{btbsystem}\rbr_{j\in\Nn}$ and the fractional system $\lbr\eqref{isstbsfrcsystem-half}\rbr_{\j\in\Nn}$ are equivalent iff the collection of functions $\lbr\sUj,\sVj\rbr_{j\in\Nn}$ satisfies \eqref{equivalence-cond}  for $ j\in\Nn$.  If $\lbr\sUj,\sVj\rbr_{j\in\Nn}$ are the BTBS functionals defined as in \thmref{A-BTBSPDE} and \thmref{AN-BTBSPDE}, then $\lbr\sUj,\sVj\rbr_{j\in\Nn}$ satisfy \eqref{equivalence-cond}  for every $j\in\Nn$.  If $n=1$, then $u=\sUj=\sVj$, the one-parameter version of the condition \eqref{equivalence-cond} is satisfied, and the equivalence between \eqref{BTBMPDE} and \eqref{BTBM-half} holds.
\end{thm}
\efr
\subsection{Inverse-stable-L\'evy-time Brownian sheets and their interacting fractional and high-order PDEs systems}\lbl{inverse}
Inverse stable subordinators---which we also call inverse L\'evy motions---arise in \cite{limitCTRW,Zsolution} as scaling limits of continuous time random walks.  Let $S(n)=Y_1+\cdots+Y_n$ a sum of independent and identically distributed random variables with $EY_n=0$ and $EY_n^2<\infty$.  The scaling limit $c^{-1/2}S([ct])\Rightarrow B(t)$ as $c\to\infty$ is a Brownian motion $B$ at time $t$, which is normal with mean zero and variance proportional to $t$.  Consider $Y_n$ to be the random jumps of a particle.  If we impose a random waiting time $T_n$ before the $n$th jump $Y_n$, then the position of the particle at time $T_n=J_1+\cdots+J_n$ is given by $S(n)$.  The number of jumps by time $t>0$ is $N(t)=\max\{n:T_n\leq t\}$, so the position of the particle at time $t>0$ is $S(N(t))$, a subordinated process.  If $\P(J_n>t)=t^{-\beta}l(t)$ for some $0<\beta<1$, where $l(t)$ is slowly varying, then the scaling limit $c^{-1/\beta}T_{[ct]}\Rightarrow L(t)$ is a strictly increasing stable L\'evy motion $L$ at time $t$ and with index $\beta$, sometimes called a stable subordinator.  The jump times $T_n$ and the number of jumps $N(t)$ are inverses $\{N(t)\geq x\}=\{T(\ceil{x})\leq t\}$ where $\ceil{x}$ is the smallest integer greater than or equal to $x$.  It follows that the scaling limits are also inverses $c^{-\beta}N(ct)\Rightarrow \Lambda(t)$ where $\Lambda(t)=\inf\{x:L(x)> t\}$, so that $\{\Lambda(t)\leq x\}=\{L(x)\geq t\}$.  We call the process $\Lambda$ a $\beta$-inverse L\'evy motion.  Since $N({ct})\approx c^{\beta}\Lambda(t)$, the particle location may, for large $c$, be approximated by  $c^{-\beta/2}S(N({[ct]}))\approx (c^\beta)^{-1/2}S(c^\beta \Lambda(t))\approx B(\Lambda(t))$, a Brownian motion subordinated to the inverse or hitting time (or first passage time) process of the stable subordinator $L$.  The random variable $L(t)$ has a smooth density.  For properly scaled waiting times, the density of $L(t)$ has Laplace transform $e^{-ts^\beta}$ for any $t>0$, and the random variables $L(t)$ and $t^{1/\beta}L(1)$ are identically distributed.  Writing $g_\beta(u)$ for the density of $L(1)$, it follows that $L(t)$ has density $t^{-1/\beta}g_\beta(t^{-1/\beta}u)$ for any $t>0$.  Using the inverse relation $\P(\Lambda(t)\leq x)=\P(L(x)\geq t)$ and taking derivatives, it follows that $\Lambda(t)$ has density
\begin{equation}\label{Etdens}
\KEmtzx=t\beta^{-1}x^{-1-1/\beta}g_\beta(tx^{-1/\beta}) ,
\end{equation}
whose $t\mapsto s$ Laplace transform $s^{\beta-1}e^{-xs^\beta}$ can also be derived from the equation
\[\KEmtzx=\frac d{dx}\P(L(x)\geq t)=\frac d{dx}\int_t^\infty x^{-1/\beta}g_\beta(x^{-1/\beta}u)\,du\]
by taking Laplace transforms on both sides.  Some fundamental properties  of $\KEmtzx$ are summarized below
\begin{lem}[Hahn et al. \cite{hahn-et-al}]\label{fractional-density-pde}
Let $\KEmtzx$ be the density of $\Lambda(t)$. Then
\begin{itemize}
\item[(a)] $\lim_{t\to+0} \KEmtzx=\delta_0(x)$ in the sense of the topology of the space of tempered distributions $\mathscr{D}'(\R)$;
    \item[(b)] $\lim _{x\to +0}\KEmtzx=\frac{t^{-\beta}}{\Gamma(1-\beta)}, \ t>0$;
    \item[(c)]$\lim _{x\to \infty}\KEmtzx=0,\ t>0$;
    \item[(d)] $t\mapsto s$ Laplace transform of $\KEmtzx$ is  $s^{\beta-1}e^{-xs^\beta}$;
    \item[(e)]For each $t>0$, $\KEmtzx$ satisfies
    \beq\label{fractional-density}
    \partial_t^\beta \KEmtzx=-\partial_x \KEmtzx-\frac{t^{-\beta}}{\Gamma(1-\beta)}\delta_0(x)
    \eeq
    in the sense of tempered distributions.
\end{itemize}
\end{lem}

Let $X$ be a Brownian motion running with twice the speed of standard Brownian motion (i.e., $X(t)=B(2t)$, where $B$ is a standard Brownian motion). The density of $X(t)$ is $\frac{1}{\sqrt{4\pi t}}e^{-|x|^2/4t}$ (we denote by $\KBmitzs=\frac{1}{\sqrt{2\pi t}}e^{-|x|^2/2t}$ the density of a one-dimensional BM starting at $0$).  From \eqref{Etdens} it follows that  for $\beta=\frac12$
\begin{equation}\begin{split}\label{Ydens}
\KEmhalftzx=2tg_{1/2}({t}/{x^{2}})x^{-3}&=\frac{2t}{x^{3} \sqrt{4\pi t^3 / x^6}}\exp\left(-\frac{x^2}{4t}\right)\\
&= \frac{2}{\sqrt{4\pi t}}\exp\left(-\frac{x^2}{4t}\right)\\
&=\frac{1}{\sqrt{2}}\KabBmtzxsqrt2
\end{split}\end{equation}
where $\KabBmtzx$ is the transition density of $|B|$.  Hence we get from the above observations that
$$
\KabBmtzx=\sqrt2\KEmhalftzxsqrt2
$$
\beq\label{fractional-density-half}
\bsp
    \partial_t^{\frac12} \KabBmtzx&=-\frac{\sqrt{2}}{2}\partial_x \KabBmtzx-\sqrt{2}\frac{t^{-\frac12}}{\Gamma(\frac12)}\delta_0(x)\\
    \pa_{t}\KEmhalftzx&=\pa^{2}_{x}\KEmhalftzx.
\end{split}
\eeq
More generally, we have the following differential facts for $\KEmtzx$ (see equation (3.18) in \cite{keyantuo-lizama})
\blm[Higher order PDEs associated with the $\Lambda$ density]\lbl{lambdapde}
For $\beta=1/\nu$, $\nu\in\lbr2,3,4,\ldots \rbr$, the density $\KEmtzx$ satisfies
\begin{equation}\label{Et-first-time-derivative}
\begin{split}
\pa_{t}\KEmtzx&= (-1)^{\nu}\pa^{\nu}_{x}\KEmtzx;\ t>0,\ x>0\\
\pa^{k}_{x}\KEmtzx\bigg|_{x=0}&=t^{-(k+1)/\nu}\frac{(-1)^k}{\Gamma (1-\frac{(k+1)}{\nu})},\ \ t>0,  k= 0,1,2,\cdots, (\nu-2); \\
\pa^{\nu -1}_{x}\KEmtzx\bigg|_{x=0}&=0,\ \ t>0;\\
\lim _{x\to\infty}\pa^{k}_{x}\KEmtzx&=0, \ \  t>0, k=0,1,2,\cdots, (\nu-1).
\end{split}
\end{equation}

\elm

Let $\Lambda^{(1)},\ldots,\Lambda^{(n)}$ be $n$ independent copies of a $\beta$-inverse L\'evy motion of index $0<\beta<1$  starting at $0$ and independent of an $n$-parameter ($n\ge1$) $\Rd$-valued Brownian sheet $W^{x}$.
We define the $n$-parameter $\Rd$-valued inverse-L\'evy-time Brownian sheet (ISLTBS) based on $\Wx$ and $\Lambda^{(1)},\ldots, \Lambda^{(n)}$, and starting at $x\in\Rd$, by
\beq\lbl{ISSTBSdef}
\WxEnt\eqdef \Wx\lpa  \Lambda^{(1)}(t_{1}),\ldots, \Lambda^{(n)}(t_{n})\rpa;\mbox{ } \t\in\Rpn.
\eeq
Clearly, $\P\lbk\WxEnt=x\rbk=1$ for $\t\in\bRpn$.  Our first result in this more general setup gives an interacting time-fractional PDEs system that is solved by running the inverse-stable-L\'evy-time Brownian sheet in \eqref{ISSTBSdef}.  It generalizes \thmref{AN-BTBSPDE} from the case $\beta=\frac12$ to that of $0<\beta<1$; and it also generalizes \thmref{A-BTBSPDE} by giving other high ($2\nu$) order interacting PDEs systems for $\beta=1/\nu$, $\nu\in\lbr2,3,4,\ldots\rbr$.
\bfr
\begin{thm}[The ISLTBS fractional and high order systems]\lbl{ISLTBSPDE} Let $0<\beta<1$.
Let $\{\WxEnt; \t\in\Rpn\}$ be an $n$-parameter $\Rd$-valued ISLTBS based on a Brownian sheet $\Wx$ and $\beta$-inverse L\'evy motions $\lbr \Lambda^{(i)}\rbr_{i=1}^{n}$ and starting on $\bRpn$ at $x\in\Rd$ on $\OFP$.
Assume the same conditions of \thmref{A-BTBSPDE} on $f:\Rd\to\R$.   If
\beq{}\label{ISSTBSexp0}
\bsp
u(\t,x)&=\EP\lbk f\lpa\WxEnt\rpa\rbk,
\\ \mathscr{V}^{(j)}(\t,x)&=\EP\lbk \lpa\prod_{i\in\Nn\setminus\lbr j\rbr} \Lambda^{(i)}(t_{i})\rpa f\lpa\WxEnt\rpa\rbk;\ j\in\Nn,
\end{split}
\eeq
for $(\t,x)\in\Rpn\times\Rd$, then the family $\lbr u, \mathscr{V}^{(j)}(\t,x)\rbr_{j\in\Nn}$ is a solution to the fractional PDEs system

 \begin{equation}\label{isstbsfrcsystem}
\begin{cases}
(a)\ \displaystyle\p_{t_j}^\beta u=\frac12\lap_{x}\mathscr{V}^{(j)};& \t\in\intRpn,x\in\Rd,\\
(b)\ u(\t,x)=f(x);&\t\in\bRpn,x\in\Rd,\\
(c)\ \mathscr{V}^{(j)}(\t,x)=0;&\tnj\in\p\Rpnmo, x\in\Rd,\\
(d)\ \mathscr{V}^{(j)}(\t,x)=\lbk \prod_{i\in\Nn\setminus\lbr j\rbr}E(\beta,1)t_{i}^\beta\rbk f(x);&\ t_{j}=0, x\in\Rd,
 \end{cases}
\end{equation}
for $ j\in\Nn$$;$
where $E(\beta,\gamma)=\EP(\Lambda(1)^\gamma)=\int_0^\infty x^{-\gamma\beta}g_\beta(x)dx$, for $\gamma>-1.$

If we now assume $\beta=1/\nu$$,$ $\nu\in\lbr2,3,4,\ldots\rbr$$,$ all the $2\nu-2$ order derivatives of $f$ are bounded and H\"older continuous, with exponent $0<\alpha\le1$, and
\beq\label{ISLTBSexp0}
\bsp
\mathscr{U}_{\nu}^{(j)}(\t,x)&=\EP\lbk \lpa\prod_{i\in\Nn\setminus\lbr j\rbr} \Lambda^{(i)}(t_{i})\rpa^{\nu}f\lpa\WxEnt\rpa\rbk;\ j\in\Nn,
\end{split}
\eeq
for $(\t,x)\in\Rpn\times\Rd$$;$ then the family $\lbr u, \mathscr{U}_{\nu}^{(j)}(\t,x)\rbr_{j\in\Nn}$ is a solution to the $2\nu$-order PDEs system
\begin{equation}\label{iltbsystem-2m}
\begin{cases}
(a)\ \displaystyle\frac{\p u}{\p{t_{j}}}=\sum_{\kappa=1}^{\nu-1}\frac{\Delta_{x}^{\kappa}f}{2^{\kappa}}\partial_{t_{j}}M_{\kappa}^{(j)}+\frac{1}{2^{\nu}}\Delta^{\nu}_{x}\sU_{\nu}^{(j)};& \t\in\intRpn,x\in\Rd,\\*
(b)\ u(\t,x)=f(x);&\t\in\bRpn,x\in\Rd,\\
(c)\ \sU_{\nu}^{(j)}(\t,x)=0;&\tnj\in\p\Rpnmo, x\in\Rd,\\
(d)\ \sU_{\nu}^{(j)}(\t,x)= f(x)N^{(j)}_{\nu}(\t) ;&\ t_{j}=0, x\in\Rd,
 \end{cases}
\end{equation}
for $ j\in\Nn$$;$ where
\beq
\bsp
N^{(j)}_{\nu}(\t)&=\prod_{i\in\Nn\setminus\lbr j\rbr}\EP\lpa\Lambda^{(i)}(t_{i})\rpa^{\nu}=E(1/\nu, \nu)^{n-1}\prod_{i\in\Nn\setminus\lbr j\rbr} t_i\\
M_{\kappa}^{(j)}(\t)& =\lbk \prod_{i\in\Nn\setminus\lbr j\rbr}\EP\lpa\Lambda^{(i)}(t_{i})\rpa^\kappa\rbk\frac{\EP\lpa\Lambda^{(j)}(t_{j})\rpa^{\kappa}}{\kappa!}\\
&=\int_{\Rpn}\lbk\prod_{i\in\Nn\setminus\lbr j\rbr}s_i\rbk^{\kappa} \frac{s_{j}^{\kappa}}{\kappa!}\prod_{i=1}^{n} \KEmitzs d\s = \frac{E(1/\nu,\kappa)^n}{\kappa !}\prod_{i\in\Nn}t_i^{ \kappa/\nu}
\end{split}
\eeq
and where  $\EP[\Lambda(t)^{\gamma}]$ for $-1<\gamma$ is given by
  $$
  \EP[\Lambda(t)^{\gamma}]=\EP[(L(1)/t)^{-\gamma/\nu}] =t^{\gamma/\nu}\EP [L(1)^{-\gamma/\nu}]=:t^{\gamma/\nu}E(1/\nu, \gamma)<\infty.
  $$
 \end{thm}
 \efr
 \subsection{Conditional equivalence of the ISLTBS high order and time-fractional interacting  PDEs systems}
Our fourth main Theorem provides the conditions needed for equivalence to hold between the high order and the fractional systems in \thmref{ISLTBSPDE}.
\bnt
The $k$-iterated $\beta$-fractional derivative in $t$ is denoted by $\partial^{k\otimes\beta}_{t}$
\ent
\bfr
 \begin{thm}[Conditional equivalence of the fractional and the $(2\nu)$-order interacting PDEs systems]\lbl{equivisltbs}
Assume that $\beta=1/\nu$$,$ $\nu\in\lbr2,3,4,\ldots\rbr$$,$ all the $2\nu-2$ order derivatives of $f$ are bounded and H\"older continuous, with exponent $0<\alpha\le1$.   Assume that, for every $j\in\Nn$, a function $u:\Rpn\times\Rd\to\R$ satisfies
 \beq\label{u-coefficients}
   \pa_{t_j}^{k\otimes\beta}u\bigg|_{t_j=0}=\frac{ \Gamma\lpa\frac{\nu-k}{\nu}\rpa E\lpa\frac1\nu,k \rpa^n\lap_x^{k}f(x)\prod_{i\in\Nn\setminus\lbr j\rbr}t_i^{ \frac k\nu}}{\nu2^{k}(k-1) !}
 \eeq
for $\ k=1,\cdots, \nu-1. $
Suppose that, for every $\nu$, the pair $\lbr\sUj_{\nu},\sVj\rbr$ satisfies
\begin{equation}\label{equivalence-cond-iltbs}
\begin{cases}
(a)\ \displaystyle2^{\nu-1}\lap_x\lpa\partial_{t_{j}}^{(\nu-1)\otimes\beta}\sVj\rpa=\lap_x^\nu\sU_{\nu}^{(j)}; & \t\in\intRpn,x\in\Rd,\\
(b)\ \sU_{\nu}^{(j)}(\t,x)= f(x)N^{(j)}_{\nu}(\t);&\ t_{j}=0, x\in\Rd,\\
(c)\ \mathscr{V}^{(j)}(\t,x)= f(x)\ds\prod_{i\in\Nn\setminus\lbr j\rbr}E\lpa\frac{1}{\nu},1\rpa t_{i}^{\beta};&\ t_{j}=0, x\in\Rd,\\
 (d)\ \mathscr{V}^{(j)}(\t,x)=\sU_{\nu}^{(j)}(\t,x)=0;&\tnj\in\p\Rpnmo, x\in\Rd,
 \end{cases}
 \end{equation}
  for some $j^{*}\in\Nn$ $($all the derivatives exist and are continuous$)$.  Then \eqref{isstbsfrcsystem} is satisfied for $j^{*}$ by $\lbr u, \mathscr{V}^{(j^{*})}\rbr$ iff \eqref{iltbsystem-2m} is satisfied for $ j^{*}$  by  $\lbr u,  \sU_{\nu}^{(j^{*})}\rbr$.   Conversely, if  \eqref{u-coefficients}, \eqref{isstbsfrcsystem}, and \eqref{iltbsystem-2m} are satisfied by $\lbr u, \mathscr{V}^{(j^{*})}\rbr$ and $\lbr u,  \sU_{\nu}^{(j^{*})}\rbr$ for $j^{*}\in\Nn$; then \eqref{equivalence-cond-iltbs} holds for $j^{*}$.  In particular, under the condition \eqref{u-coefficients}, the high order system $\lbr\eqref{iltbsystem-2m}\rbr_{j\in\Nn}$ and the fractional system $\lbr\eqref{isstbsfrcsystem}\rbr_{\j\in\Nn}$ are equivalent iff the collection of functions $\lbr\sUj_{\nu},\sVj\rbr_{j\in\Nn}$ satisfies \eqref{equivalence-cond-iltbs} for $ j\in\Nn$.  If $\lbr\sUj_{\nu},\sVj\rbr_{j\in\Nn}$ are the ISLTBS functionals defined as in \thmref{ISLTBSPDE}, then $\lbr\sUj_{\nu},\sVj\rbr_{j\in\Nn}$ satisfy \eqref{equivalence-cond-iltbs} and \eqref{u-coefficients} for every $j\in\Nn$.  If $n=1$, then $u=\sUj=\sVj$, the one-parameter version of the conditions \eqref{equivalence-cond-iltbs} and \eqref{u-coefficients} are satisfied and the equivalence between \eqref{isstbsfrcsystem} and \eqref{iltbsystem-2m} holds.
\end{thm}
\efr
\subsection{Examples of the initial function $f$}\lbl{fexample}
Theorem 4.5 in Gilbarg and Trudinger \cite{gilbarg-trudinger} gives the fact that when $f_{0}$ is H\"older continuous of order $0<\alpha\le1$ with compact support and $f\in C_{c}^2(\R^d)$  (with compact support) satisfy Poisson's equation $\Delta f=f_{0}$, then $f\in C_{c}^{2,\alpha}(\R^d)$ (the space of functions such that all the derivatives of order  $2$ are H\"older continuous of order $0<\alpha\le1$  with compact support).
Similarly, if $f_{0}\in C_{c}^{k, \alpha}(\R^d)$ and $f\in C_{c}^2(\R^d)$ satisfy $\Delta f=f_{0}$, then $f\in C_{c}^{k+2, \alpha}(\R^d)$.  The condition of compact support can be removed by various means (see Theorem 4.6 in Gilbarg and Trudinger \cite{gilbarg-trudinger} and the discussion before).

We now give two related types of functions satisfying the conditions on the initial function $f$ in our results.  Let $C$ be an arbitrary constant. Let $f_{0}$ :$\Rd\to\R$ be given by
\beq\lbl{fexmpl}
f_{0}(x)=\bc C\exp\lpa\df{1}{\lab x\rab^{2\alpha}-1}\rpa;&\lab x\rab<1,\cr
0;&|x|\ge1,
\ec
\eeq
then $f_0(x)$ is $\alpha$-H\"older continuous with compact support.  When $\alpha=1$, $f_{0}\in C^{\infty}_{c}\lpa\Rd\rpa$ with all the derivatives being Lipschitz.  We can then set $f=f_{0}$ in \thmref{A-BTBSPDE}, \thmref{AN-BTBSPDE},  and  \thmref{ISLTBSPDE} to get the corresponding PDEs.  In the case of \thmref{A-BTBSPDE}, if $d=2$, then $$\Delta_{x}f(x)=\df{4C\lpa|x|^{4}+|x|^{2}-1\rpa\exp\lpa\df{1}{\lab x\rab^{2}-1}\rpa}{\lpa|x|^{2}-1\rpa^{4}}$$
Alternatively, iterating the Gilbarg and Trudinger procedure above, we can use $f_0$ in \eqref{fexmpl} to generate functions whose derivatives of order $2\nu-2 $  are all $\alpha$-H\"older continuous for $\nu=\frac1\beta\in\lbr2,3,\ldots\rbr$.  In particular, we can pick $f$ such that $\Delta f=f_0$  and $f\in C_{c}^{2,\alpha}(\R^d)$. In this case, with \beq\label{BTBSexp0-particular}
\bsp
u(\t,x)&=\EP\lbk f\lpa\WxBnt\rpa\rbk,
\\ \mathscr{U}^{(j)}(\t,x)&=\EP\lbk \lpa\prod_{i\in\Nn\setminus\lbr j\rbr}\lab B^{(i)}(t_{i})\rab\rpa^{2}f\lpa\WxBnt\rpa\rbk;\ j\in\Nn,
\end{split}
\eeq
for $(\t,x)\in\Rpn\times\Rd$, we have that the family $\lbr u,  \sU^{(j)}\rbr_{j\in\Nn}$ is a solution to 
the system of interacting fourth order linear PDEs:
\begin{equation}\label{btbsystem-particular}
\begin{cases}
(a)\ \displaystyle{\p_{t_{j}} u}=\sqrt{\frac{\prod_{i\in\Nn\setminus\lbr j\rbr}t_i}{2^{4-n}t_j\pi^{n}}}f_0(x)+\frac18\lap_{x}^{2}\sU^{(j)};& \t\in\intRpn,x\in\Rd,\\
(b)\ u(\t,x)=f(x);&\t\in\bRpn,x\in\Rd,\\
(c)\ \sU^{(j)}(\t,x)=0;&\tnj\in\p\Rpnmo, x\in\Rd,\\
(d)\ \sU^{(j)}(\t,x)=\lbk \prod_{i\in\Nn\setminus\lbr j\rbr}t_{i}\rbk f(x);&\ t_{j}=0, x\in\Rd,
 \end{cases}
\end{equation}
for $ j\in\Nn.$
We can also write the results in \thmref{ISLTBSPDE} in a similar way.

\section{Establishing the BTBS fourth-order and fractional PDEs systems links}
\subsection{A stochastic analytic unifying proof of \thmref{AN-BTBSPDE} and \thmref{A-BTBSPDE}}\lbl{newsapf}
First, we need the following lemma that allows us to pull the Laplacian and bi-Laplacian in and out of different integrals.
\bnt\lbl{laplacians}  We bring the reader's attention to the difference between the spatial Laplacian
\beq\lbl{slap}
\Delta_{x}\phi(\t,x)=\sum_{k=1} ^{d}\pd^{2}_{x_{k}}\phi(\t,x)
\eeq
and the Laplacian
\beq\lbl{Itolap}
\Delta g\lpa w(\t,x)\rpa=\sum_{k=1} ^{d}\pd^{2}_{y_{k}}g\lpa y\rpa|_{y_{k}=w_{k}(\t,x)},
\eeq
 for sufficiently smooth $\phi:\Rpn\times\Rd\to\R$ and $g:\Rd\to\R$, and for $w:\Rpn\times\Rd\to\Rd$.  The second Laplacian \eqref{Itolap} appears when we apply It\^o's formula in our proofs.
\ent
\blm[Interchanging $\Delta$, $\Delta_{x}$, and $\Delta_{x}^{2}$ with different integrals: the BTBS case]\lbl{lapintBTBS}
Assume the same conditions on $f$ as in \thmref{A-BTBSPDE}.  Then, recalling the notational conventions in \notnref{papernot}, we have
\ben\rencomalp
\item $\EP\lbk\Delta f\lpa\Wx_{\snj}(\rho_{j})\rpa\rbk =\Delta_{x}\EP\lbk f\lpa\Wx_{\snj}(\rho_{j})\rpa\rbk,$
\item $\ds\Delta_{x}\int_{0}^{r_{j}}\EP\lbk\Delta f\lpa\Wx_{\snj}(\rho_{j})\rpa\rbk d\rho_{j}=\int_{0}^{r_{j}}\Delta^{2}_{x}\EP\lbk f\lpa\Wx_{\snj}(\rho_{j})\rpa\rbk d\rho_{j},$
\item
\beqs
\bsp
&\lap_{x}\int_{\Rpn}\lbk\prod_{i\in\Nn\setminus\lbr j\rbr}s_i\rbk \EP\lbk f\lpa\Wx_{\s}\rpa\rbk\prod_{i=1}^{n}\KBmitzs d\s
\\=&\int_{\Rpn}\lbk\prod_{i\in\Nn\setminus\lbr j\rbr}s_i\rbk \lap_{x}\EP\lbk f\lpa\Wx_{\s}\rpa\rbk\prod_{i=1}^{n}\KBmitzs d\s,
\end{split}
\eeqs
and
\item \beqs
\bsp
&\lap_{x}^{2}\int_{\Rpn}\lbk\prod_{i\in\Nn\setminus\lbr j\rbr}s_i^{2}\rbk \EP\lbk f\lpa\Wx_{\s}\rpa\rbk\prod_{i=1}^{n}\KBmitzs d\s
\\=&\int_{\Rpn}\lbk\prod_{i\in\Nn\setminus\lbr j\rbr}s_i^{2}\rbk \lap_{x}^{2}\EP\lbk f\lpa\Wx_{\s}\rpa\rbk\prod_{i=1}^{n}\KBmitzs d\s.
\end{split}
\eeqs
\een
\elm
\brm\lbl{pulrout}
In parts (a) and (b) we wrote $\Wx_{\snj}(\rho_{j})$ in place of $\Wx_{\s}$ since this is the form that appears in the main proof of \thmref{A-BTBSPDE} and \thmref{AN-BTBSPDE} upon applying It\^o's formula in \eqref{BTBSreveal}, making it more convenient for the reader to see the correspondence.
\erm
 \bpf  Let $\KBsnjrjxy$ be the density of a Brownian sheet $\Wx_{\snj}(r_{j})$.  Part (a) follows by integration by parts twice and the facts that
\beqs
\lim_{y_{k}\to\pm\infty}\pd_{y_{k}}f(y)\KBsnjrjxy=\lim_{y_{k}\to\pm\infty}f(y)\pd_{y_{k}}\KBsnjrjxy=0
\eeqs
which give
\beq\lbl{outexpec}
\bsp
\EP\lbk\Delta f(\Wx_{\snj}(r_{j}))\rbk&=\int_{\Rd}\sum_{k=1} ^{d}\pd^{2}_{y_{k}}f(y)\KBsnjrjxy dy
\\&=\sum_{k=1}^{d}\int_{\Rd}f(y)\pd^{2}_{x_{k}}\KBsnjrjxy dy
\\&=\Delta_{x}\EP\lbk f(\Wx_{\snj}(r_{j}))\rbk,
\end{split}
\eeq
where the last equality follows from standard analysis and the fact that $$\int_{\Rd}\lab f(y)\pd^{2}_{x_{k}}\KBsnjrjxy \rab dy<\infty$$ for each $k=1,\ldots,d$ by the boundedness of $f$.  Part (b) is easily seen by noticing that, under the boundedness of the second partial derivatives of $f$, $h_{\snj}(\rho_{j},x)=\EP\lbk\Delta f\lpa\Wx_{\snj}(\rho_{j})\rpa\rbk$ is twice continuously differentiable in $x$; in fact, by the Brownian sheet-PDE connection in Lemma 1.1 in \cite{Abtbs} $h_{\snj}$ is the unique bounded classical solution to the initial value heat equation with parametrized conduction coefficient
\beq\lbl{bspde}
\bc
\ds\pa_{\rho_{j}}h_{\snj}=\frac12\lbk\prod_{i\in\Nn\setminus\{j\}}s_{i}\rbk\lap_{x}h_{\snj};&\rho_{j}>0, x\in\Rd,\\
h_{\snj}(0,x)=\lap_{x}f(x);&x\in\Rd.
\ec
\eeq
Thus, $\ds\int_{0}^{r_{j}}\lab\Delta_{x}\EP\lbk\Delta f\lpa\Wx_{\snj}(\rho_{j})\rpa\rbk\rab d\rho_{j}<\infty$, and hence
$$\ds\Delta_{x}\int_{0}^{r_{j}}\EP\lbk\Delta f\lpa\Wx_{\snj}(\rho_{j})\rpa\rbk d\rho_{j}=\int_{0}^{r_{j}}\Delta_{x}\EP\lbk\lap f\lpa\Wx_{\snj}(\rho_{j})\rpa\rbk d\rho_{j},$$
which when combined with part (a) implies (b).
Part (d) is exactly Lemma 2.1 in \cite{Abtbs} and part (c) follows the same steps with obvious modifications (see the more general case in \lemref{lapintISLTBS} (a)).
 \epf

We now give a proof that simultaneously establishes \thmref{A-BTBSPDE} and \thmref{AN-BTBSPDE}, capturing the fourth order PDE \eqref{btbsystem} and the fractional PDE \eqref{isstbsfrcsystem-half} in one equation.
\bpfs{Proof of \thmref{A-BTBSPDE} and \thmref{AN-BTBSPDE}}  From our definition of the BTBS $\WxBn$ we have that $B^{(1)},\ldots,B^{(n)}$ are $n$ independent copies of a standard one-dimensional Brownian motion starting at $0$, and we also have that $B^{(1)},\ldots,B^{(n)}$ are independent of the $n$-parameter $\Rd$-valued Brownian sheet $\Wx$.  Using these facts, we easily get that 
\beq\lbl{BTBSexplicit}
\bsp
&\mathrm{(a)}\ u(\t,x)=\EP\lbk f\lpa\WxBnt\rpa\rbk
=2^{n}\int_{\Rpn}\EP\lbk f(\Wx_{\s})\rbk\prod_{i=1}^{n}\KBmitzs d\s,\\
&\mathrm{(b) }\  \sU^{(j)}(\t,x)=\EP\lbk \lpa\prod_{i\in\Nn\setminus\lbr j\rbr}\lab B^{(i)}(t_{i})\rab\rpa^{2}f\lpa\WxBnt\rpa\rbk\\
&\hspace{2.1cm}=2^{n}\int_{\Rpn}\lbk\prod_{i\in\Nn\setminus\lbr j\rbr}s_i^{2}\rbk\EP\lbk f\lpa\Wx_{\s}\rpa\rbk \prod_{i=1}^{n}\KBmitzs d\s, \\
&\mathrm{(c) }\  \sVj(\t,x)=\EP\lbk \lpa\prod_{i\in\Nn\setminus\lbr j\rbr}\lab B^{(i)}(t_{i})\rab\rpa f\lpa\WxBnt\rpa\rbk\\
&\hspace{2.0cm}=2^{n}\int_{\Rpn}\lbk\prod_{i\in\Nn\setminus\lbr j\rbr}s_i\rbk\EP\lbk f\lpa\Wx_{\s}\rpa\rbk \prod_{i=1}^{n}\KBmitzs d\s, 
\end{split}
\eeq
for $j\in\Nn$.   Now, fix an arbitrary $j\in\Nn$, and let
\beq\lbl{vsnjrj}
v_{\snj}(\rho_{j},x):=\EP\lbk f(\Wx_{\snj}(\rho_{j}))\rbk ;\ 0\le \rho_{j}\le s_{j}.
\eeq
We apply It\^o's formula twice to $f\lpa W_{\s}^{x}\rpa$ in \eqref{BTBSexplicit} (a) in the ${j}$-th temporal parameter (see the detailed discussion in Appendix A in \cite{Abtbs}); using the independence of the coordinates of the Brownian-sheet $\Wx$, the fact that $\lbr\Wx_{\snj}(\rho_{j}); 0\le \rho_{j}<\infty\rbr$ is a martingale with quadratic variation $\lqv\Wx_{\snj}(\cdot)\rqv_{r_{j}}=r_{j}\prod_{i\in\Nn\setminus\lbr j\rbr}s_i$, and taking expectations and using the fact that the expectation of the stochastic integrals are zero,  along with \lemref{lapintBTBS} (a) and (b), we get

\beq\label{BTBSreveal}
\bsp
&u(\t,x)
=2^{n}\int_{\Rpn}\EP\lbr f\lpa\Wsnjzx\rpa+\int_{0}^{s_{j}}\nabla f\lpa\Wsnjrx\rpa\cdot d\Wsnjrx\right.\\& \left.\hspace{3cm}+\df12\lbk\prod_{i\in\Nn\setminus\lbr j\rbr}s_i\rbk\int_{0}^{s_{j}}\Delta f\lpa\Wsnjrx\rpa dr_{j}\rbr \prod_{i=1}^{n}\KBmitzs d\s
\\&=f(x)+\frac{1}{2}\lbr2^{n}\int_{\Rpn}\lbk\prod_{i\in\Nn\setminus\lbr j\rbr}s_i\rbk \int_{0}^{s_{j}}\Delta_{x}v_{\snj}(r_{j},x)dr_{j}\prod_{i=1}^{n}\KBmitzs d\s\rbr
\\&=f(x)+\frac{1}{2}\lbr2^{n}\int_{\Rpn}\lbk\prod_{i\in\Nn\setminus\lbr j\rbr}s_i\rbk \int_{0}^{s_{j}}\Delta_{x}\EP\lbr f\lpa\Wsnjzx\rpa\right.\right.\\&\left.\left.\hspace{2.1cm}+\int_{0}^{r_{j}}\nabla f\lpa\Wx_{\snj}(\rho_{j})\rpa\cdot d\Wx_{\snj}(\rho_{j})\right.\right.\\& \left.\left.\hspace{2.1cm}+\df12\lbk\prod_{i\in\Nn\setminus\lbr j\rbr}s_i\rbk\int_{0}^{r_{j}}\Delta f\lpa\Wx_{\snj}(\rho_{j})\rpa d\rho_{j}\rbr dr_{j}\prod_{i=1}^{n}\KBmitzs d\s\rbr
\\&=f(x)+\sqrt{\frac{2^{n-2}\prod_{i=1}^{n}t_{i}}{\pi^{n}}}\Delta_{x}f(x)
\\&+\frac{1}{4}\lbr 2^{n}\int_{\Rpn}\lbk\prod_{i\in\Nn\setminus\lbr j\rbr}s_i\rbk^{2}\int_{0}^{s_{j}}\int_{0}^{r_{j}}\Delta_{x}^{2}v_{\snj}(\rho_{j},x)d\rho_{j}dr_{j}\prod_{i=1}^{n}\KBmitzs d\s\rbr
\end{split}
\eeq

Equation \eqref{BTBSreveal} contains in it the two BTBS-PDEs facts that $\lbr u,\sUj\rbr$ solves the fourth-order PDE in \eqref{btbsystem} and $\lbr u,\sVj\rbr$ solves the half-derivative PDE in \eqref{isstbsfrcsystem-half}.  We now show this fact.

The interacting fourth order PDEs system in \thmref{A-BTBSPDE} is seen by first taking the $j$-th time derivative $\partial_{t_{j}}$, with respect to $t_{j}$, using the last equality in \eqref{BTBSreveal} and the dominated convergence theorem to bring $\partial_{t_{j}}$ under the integral.  We then use the fact that $\partial_{t_{j}}\KBmjtzs=\frac12\partial_{s_{j}}^{2}\KBmjtzs$ together with two integration by parts in which the boundary terms vanish.  This last fact is because of the exponential decay of $\partial_{s_{j}}\KBmjtzs$ and $\KBmjtzs$ which nullify the polynomial growth of the inside integrals that is due to the continuity of $\Delta_{x}^{2}v_{\snj}(\rho_{j},x)=\Delta_{x}\EP\lbk\Delta f\lpa\Wx_{\snj}(\rho_{j})\rpa\rbk$---which follows from \lemref{lapintBTBS} (a), the boundedness of $\Delta f$, and the discussion leading to the Brownian sheet PDE \eqref{bspde}---and hence its boundedness on $[0,r_{j}]$.  We finish with two applications of the fundamental theorem of calculus, along with \lemref{lapintBTBS} (d), to get
\beq\lbl{fourthBTBS}
\bsp
&\partial_{t_{j}}u(\t,x)=\sqrt{\frac{\prod_{i\in\Nn\setminus\{j\}}t_{i}}{2^{4-n}t_{j}\pi^{n}}}\Delta_{x}f(x)
\\&+\frac{1}{4}\lbr 2^{n}\int_{\Rpn}\lbk\prod_{i\in\Nn\setminus\lbr j\rbr}s_i^{2}\rbk\int_{0}^{s_{j}}\int_{0}^{r_{j}}\Delta_{x}^{2}v_{\snj}(\rho_{j},x)d\rho_{j}dr_{j} \partial_{t_{j}}\prod_{i=1}^{n}\KBmitzs d\s\rbr
\\&=\sqrt{\frac{\prod_{i\in\Nn\setminus\{j\}}t_{i}}{2^{4-n}t_{j}\pi^{n}}}\Delta_{x}f(x)
\\&+\frac{\Delta_{x}^{2}}{8}\lbr 2^{n}\int_{\Rpn}\lbk\prod_{i\in\Nn\setminus\lbr j\rbr}s_i^{2}\rbk\EP\lbk f\lpa\Wx_{\s}\rpa\rbk \prod_{i=1}^{n}\KBmitzs d\s\rbr
\\&=\sqrt{\frac{\prod_{i\in\Nn\setminus\lbr j\rbr}t_i}{2^{4-n}t_j\pi^{n}}}\lap_{x} f(x)+\frac18\lap_{x}^{2}\sU^{(j)}(\t,x).
\end{split}
\eeq
On the other hand, to see the interacting fractional PDEs system of \thmref{AN-BTBSPDE}, we take the $j$-th time Caputo  half-derivative $\partial^{\frac12}_{t_{j}}$, with respect to $t_{j}$, using the next to last equality in \eqref{BTBSreveal} and the dominated convergence theorem to bring $\partial^{\frac12}_{t_{j}}$ under the integral; we then use the fact that $\partial_{t_{j}}^{\frac12}\KBmjtzs=-\frac{1}{\sqrt{2}}\partial_{s_{j}}\KBmjtzs-\frac{\sqrt{2}}{\sqrt{\pi
t_{j}}}\delta_{0}(s_{j})$\footnote{$\delta_{0}(s_{j})=\lim_{tj\searrow0}\KBmjtzs$
in the sense of tempered distributions ($\delta_{0}(s_{j})ds_{j}$ is the
Dirac measure at $0$).} (this was proved in equation \eqref{fractional-density-half} and the discussion preceding it) together with integration by parts, observing that $\int_{0}^{\infty}\int_{0}^{s_{j}}\Delta_{x}v_{\snj}(\rho_{j},x)d\rho_{j}\delta_{0}(s_{j})ds_{j}=0$, along with \lemref{lapintBTBS} (c), to get\beq\lbl{halfderBTBS}
\bsp
\partial_{t_{j}}^{\frac12}u(\t,x)&=\frac{1}{2}\lbr 2^{n}\int_{\Rpn}\lbk\prod_{i\in\Nn\setminus\lbr j\rbr}s_i\rbk \int_{0}^{s_{j}}\Delta_{x}v_{\snj}(\rho_{j},x)d\rho_{j}\partial^{\frac12}_{t_{j}} \prod_{i=1}^{n}\KBmitzs d\s\rbr
\\&=\frac{\Delta_{x}}{\sqrt{8}}\lbr 2^{n}\int_{\Rpn}\lbk\prod_{i\in\Nn\setminus\lbr j\rbr}s_i\rbk \EP\lbk f\lpa\Wx_{\s}\rpa\rbk \prod_{i=1}^{n}\KBmitzs d\s \rbr
\\&=\frac{1}{\sqrt{8}}\Delta_{x}\sVj(\t,x)
\end{split}
\eeq
The temporal-boundary conditions in \eqref{btbsystem} and \eqref{isstbsfrcsystem-half} follow easily from the definitions of $\lbr u,\sUj,\sVj\rbr$ and the proof is complete since $j$ is arbitrary.
\epfs
\subsection{Proof of the conditional equivalence result: the BTBS case}
We now establish the conditional equivalence result for the BTBS fourth order and fractional interacting systems.
\bpfs{Proof of \thmref{equiv}}
Fix an arbitrary $j\in\Nn$ and suppose that two pairs $\lbr u,\sVj\rbr$ and $\lbr u,\sUj\rbr$ satisfy \eqref{btbsystem}.and \eqref{isstbsfrcsystem-half}. Take the Fourier-Laplace transform ($x\to k$, and $t_j\to s_j$) of \eqref{isstbsfrcsystem-half} to get
 \begin{equation}\label{l-f-trans-1}
 s^{1/2}_j\bar{u}_{\tnj}(s_j,k)-s^{-1/2}_j\hat{f}(k)=\frac{-|k|^2}{\sqrt8}{\bar{\mathscr{V}}_{\tnj}^{(j)}}(s_j,k).
 \end{equation}
 Now,  multiply both sides of equation \eqref{l-f-trans-1} by $s_j^{1/2}$ to get
\begin{equation}\label{l-f-trans-2}
 s_j\bar{u}_{\tnj}(s_j,k)-\hat{f}(k)=\frac{-s_j^{1/2}|k|^2}{\sqrt8}{\bar {\mathscr{V}}_{\tnj}^{(j)}}(s_j,k).
 \end{equation}
Adding and subtracting $s_j^{-1/2}\sqrt{\frac{\prod_{i\in\Nn\setminus\lbr j\rbr}t_i}{2^{4-n}\pi^{n}}}\Gamma(1/2)(-|k|^2\hat f(k))$ to
\eqref{l-f-trans-2} gives
\begin{equation}\label{l-f-trans-3}\begin{split}
 s_j\bar{u}_{\tnj}(s_j,k)-\hat{f}(k)&=\frac{-s_j^{1/2}|k|^2}{\sqrt8} {\bar{\mathscr{V}}_{\tnj}^{(j)}}(s_j,k) \\ &-s_j^{-1/2}\sqrt{\frac{\prod_{i\in\Nn\setminus\lbr j\rbr}t_i}{2^{4-n}\pi^{n}}}\Gamma(1/2)(-|k|^2\bar f(k))\\
 &+s_j^{-1/2}\sqrt{\frac{\prod_{i\in\Nn\setminus\lbr j\rbr}t_i}{2^{4-n}\pi^{n}}}\Gamma(1/2)(-|k|^2\bar f(k)).
 \end{split}
 \end{equation}
 Observing that
 \beq
 \bsp
 \lim_{t_j\to 0}\lap_x\mathscr{V}^{(j)}(\t,x)&=\lbk \prod_{i\in\Nn\setminus\lbr j\rbr}\frac{\sqrt{2}}{\sqrt{\pi}}t_{i}^{1/2}\rbk \lap_xf(x)\\&=\sqrt{\frac{2^{n-1}\prod_{i\in\Nn\setminus\lbr j\rbr}t_i}{\pi^{n}}}\Gamma(1/2) \lap_xf(x); x\in\Rd,
 \end{split}
 \eeq
equation \eqref{l-f-trans-3} inverts to
\begin{equation}\label{new-fractional-time-heat}
 \displaystyle\frac{\p u}{\p t_j}=\sqrt{\frac{\prod_{i\in\Nn\setminus\lbr j\rbr}t_i}{2^{4-n}t_j\pi^{n}}} \lap_x f(x)+\frac{1}{\sqrt8}\lap_x(\partial_{t_j}^{1/2}\mathscr{V}^{(j)}(\t,x)); \t\in\intRpn,x\in\Rd.
 \end{equation}
Comparing equation \eqref{new-fractional-time-heat} and equation \eqref{btbsystem} (a) we must have
 $$
 \frac{1}{\sqrt8}\lap_x(\partial_{t_j}^{1/2}\mathscr{V}^{(j)}(\t,x))=\frac18\lap_x^2\mathscr{U}^{(j)}(\t,x) ; \t\in\intRpn,x\in\Rd
 $$
and condition \eqref{equivalence-cond} is satisfied.

For the rest of the proof, suppose that \eqref{equivalence-cond} holds.  We then assume first that $\lbr u,\sVj\rbr$ is a solution to \eqref{isstbsfrcsystem-half}.
   Observe that
 $$\partial_{t_j}^{1/2}(\partial_{t_j}^{1/2}u(\t,x))=\partial_{t_j}u(\t,x)-\sqrt{\frac{\prod_{i\in\Nn\setminus\lbr j\rbr}t_i}{2^{4-n}t_j\pi^{n}}} \lap_x f(x)$$
This follows from the fact that, under the smoothness conditions on $u$, it satisfies
$$
\partial_{t_j}^{1/2}(\partial_{t_j}^{1/2}u(\t,x))=\partial_{t_j}u(\t,x)-\frac{t_j^{-1/2}}{\sqrt{\pi}}\left(\partial_{t_j}^{1/2}u(\t,x)|_{t_j=0}\right)
$$
and the fact that
$$
\partial_{t_j}^{1/2}u(\t,x)|_{t_j=0}=\frac{{1}}{\sqrt8}\lim_{t_j\to 0}\lap_x \mathscr{V}^{(j)}(\t,x)=\frac{{1}}{\sqrt8}\sqrt{\frac{\prod_{i\in\Nn\setminus\lbr j\rbr}t_i}{2^{1-n}
\pi^{n-1}}} \lap_x f(x).
$$
Taking the half derivative with respect to $t_j$ of the equation in \eqref{isstbsfrcsystem-half}
 then we get that $\lbr u,\sUj\rbr$ is a solution to \eqref{btbsystem}.

 If we now assume that $\lbr u,\sUj\rbr$ is a solution to \eqref{btbsystem}, then we can use condition \eqref{equivalence-cond} (a) to rewrite \eqref{btbsystem} as
\beq \label{u-v-pde}
 \displaystyle\frac{\p u}{\p t_j}=\sqrt{\frac{\prod_{i\in\Nn\setminus\lbr j\rbr}t_i}{2^{4-n}t_j\pi^{n}}}\lap_{x} f(x)+\frac{1}{\sqrt{8}}\lap_x(\partial_{t_j}^{1/2}\mathscr{V}^{(j)}(\t,x)); \t\in\intRpn,x\in\Rd
 \eeq
  Taking the Fourier-Laplace transform ($x\to k$, and $t_j\to s_j$) of \eqref{u-v-pde} using \eqref{equivalence-cond} (d)  we get

 \beq \label{FL-transform}\begin{split}
 s_j\bar{u}_{\tnj}(s_j,k)-\hat{f}(k)&=  \sqrt{\frac{\prod_{i\in\Nn\setminus\lbr j\rbr}t_i}{2^{4-n}\pi^{n}}}s^{-1/2}_j \Gamma(1/2)(-|k|^2\bar{f}(k))\\
 &+\frac{1}{\sqrt{8}} \lbr -|k|^2s_j^{1/2}
 {\bar{\mathscr{V}}_{\tnj}^{(j)}}(s_j,k)-s_j^{-1/2}(-|k|^2)\hat{f}(k)\sqrt{\frac{\prod_{i\in\Nn\setminus\lbr j\rbr}t_i}{2^{1-n}
\pi^{n-1}}}\rbr.
\end{split}
\eeq
 Observe that the first and third terms  on the right hand side of \eqref{FL-transform} cancel each other.  Hence \eqref{FL-transform} reduces to
 \beq
  s_j\bar{u}_{\tnj}(s_j,k)-\hat{f}(k)= \frac{1}{\sqrt{8}} ( -|k|^2)s_j^{1/2} {\bar{\mathscr{V}}_{\tnj}^{(j)}}(s_j,k).
 \eeq
 Dividing both sides by $s_j^{1/2}$ we obtain
 $$
 s_j^{1/2}\bar{u}_{\tnj}(s_j,k)-s_j^{-1/2}\hat{f}(k)= \frac{1}{\sqrt{8}} ( -|k|^2) {\bar{\mathscr{V}}_{\tnj}^{(j)}}(s_j,k),
 $$
 which, upon inversion, implies that $\lbr u,\sVj\rbr$ is a solution to \eqref{isstbsfrcsystem-half}.
 \epfs
 \section{Establishing the ISLTBS high-order and fractional PDEs systems links}
\subsection{A stochastic analytic proof of \thmref{ISLTBSPDE}}\lbl{newsapf}
Like the BTBS case, we need a lemma that allows us to pull the Laplacian and bi-Laplacian in and out of different integrals.  Let $\Lambda(t)$ be the inverse of $L(t)$ with index $0<\beta<1$.
 We can calculate  $\EP[\Lambda(t)^{\gamma}]$ for $-1<\gamma$ as follows:
   $\Lambda(t)\stackrel{(d)}{=} (L(1)/t)^{-1/\nu}$ by Corollary 3.1 in Meerschaert and Scheffler \cite{limitCTRW} hence
  $$
  \EP[\Lambda(t)^{\gamma}]=\EP[(L(1)/t)^{-\gamma\beta}] =t^{\gamma\beta}\EP[L(1)^{-\gamma\beta}]=:t^{\gamma\beta}E(\beta, \gamma)<\infty.
  $$

 \blm[Interchanging $\Delta_{x}$ and $\Delta_{x}^{\nu}$ with different integrals: the ISLTBS case]\lbl{lapintISLTBS}
Let $f:\Rd\to\R$ be bounded and measurable such that
\ben\rencomalp
\item  $f$ bounded and H\"older continuous, with exponent $0<\alpha\le1$, then we have
\beqs
\bsp
&\lap_{x}\int_{\Rpn}\lbk\prod_{i\in\Nn\setminus\lbr j\rbr}s_i\rbk \EP\lbk f\lpa\Wx_{\s}\rpa\rbk\prod_{i=1}^{n} \KEmitzs d\s
\\=&\int_{\Rpn}\lbk\prod_{i\in\Nn\setminus\lbr j\rbr}s_i\rbk \lap_{x}\EP\lbk f\lpa\Wx_{\s}\rpa\rbk\prod_{i=1}^{n} \KEmitzs d\s;
\end{split}
\eeqs
\item all the $2\nu-2$ order derivatives of $f$ are bounded and H\"older continuous, with exponent $0<\alpha\le1$, then we have
\beqs
\bsp
&\Delta^{\nu}_{x}\int_{\Rpn}\lbk\prod_{i\in\Nn\setminus\lbr j\rbr}s_i^{\nu}\rbk \EP\lbk f\lpa\Wx_{\s}\rpa\rbk\prod_{i=1}^{n} \KEmitzs d\s
\\=&\int_{\Rpn}\lbk\prod_{i\in\Nn\setminus\lbr j\rbr}s_i^{\nu}\rbk \Delta^{\nu}_{x}\EP\lbk f\lpa\Wx_{\s}\rpa\rbk\prod_{i=1}^{n} \KEmitzs d\s.
\end{split}
\eeqs
\een
\elm
\bpf
The proof of (a) is similar to the proof of Lemma 2.1 in \cite{Abtbs}. We need to consider second derivative instead of fourth derivative and the function $f$ instead of the second derivative of $h$ there. We have
\begin{equation}\label{second-order-derivative-BS}
\begin{split}
\partial^2_{x_k}\EP\lbk f\lpa\Wx_{\s}\rpa\rbk & =\int_{\Rd} f(y)\partial^2_{x_k} K^{BS(n,d)}_{s;x,y}dy\\
&=\EP\lbk \bigg( \frac{(x_k-W_k^{x_k}(\s))^2-\prod_{i=1}^n s_i}{\prod_{i=1}^n s_i^2}\bigg)(f\lpa\Wx_{\s}\rpa-f(x))\rbk.
\end{split}
\end{equation}
where $K^{BS(n,d)}_{s;x,y}$ is the BS density previously given in \eqref{bsdensity}.  By scaling we obtain
$$
\EP \bigg| (x_k-W_k^{x_k}(\s))^2-\prod_{i=1}^n s_i\bigg|^2=\prod_{i=1}^n s_i^2 \EP \bigg| \bigg(\frac{W_k^{0}(\s)}{\sqrt{\prod_{i=1}^n s_i}}\bigg)^2-1\bigg|^2=C\prod_{i=1}^n s_i^2.
$$
for some constant $C$.
Using \eqref{second-order-derivative-BS} and Cauchy-Schwarz inequality, along with the fact that $f$ is H\"older continuous with exponent $\alpha$, we get
\begin{equation}
\bigg|\partial^2_{x_k}\EP\lbk f\lpa\Wx_{\s}\rpa\rbk\bigg|\leq \frac{C}{\prod_{i=1}^n s_i^{1-\alpha/2}}
\end{equation}

 Now,
\begin{equation}
\begin{split}
&\int_{\Rpn}\lbk\prod_{i\in\Nn\setminus\lbr j\rbr}s_i\rbk \bigg|\partial^2_{x_k}\EP\lbk f\lpa\Wx_{\s}\rpa\rbk\bigg|\prod_{i=1}^{n}
\KEmitzs d\s\\
&\leq \int_{\Rpn}\lbk\prod_{i\in\Nn\setminus\lbr j\rbr}s_i\rbk \frac{C}{\prod_{i=1}^n s_i^{1-\alpha/2}}\prod_{i=1}^{n}
\KEmitzs d\s\\
&=C\int_{\Rpn}\lbk\prod_{i\in\Nn\setminus\lbr j\rbr}s_i^{\alpha/2}\rbk s_j^{-(1-\alpha/2)} \prod_{i=1}^{n}
\KEmitzs d\s\\
&=C\bigg(\int_{0}^\infty s_j^{\alpha/2-1}\KEmjtzs ds_{j}\bigg) \prod_{i\in\Nn\setminus\lbr j\rbr}\int_{0}^\infty s_i^{\alpha/2}\KEmitzs d s_{i}\\
&=C \EP[\Lambda(t_j)^{\alpha/2-1}]\prod_{i\in\Nn\setminus\lbr j\rbr}\EP[\Lambda(t_i)^{(\alpha/2)}]\\
&=C t_j^{(\alpha/2-1)/\nu}E(1/\nu, \alpha/2-1)\prod_{i\in\Nn\setminus\lbr j\rbr} t_i^{\alpha/(2\nu)}E(1/\nu, \alpha/2)<\infty
\end{split}
\end{equation}

  We next prove part (b). For notational simplicity we show that
  \beq\label{k-derivative-nu-order}
\bsp
&\partial^{2\nu}_{x_k}\int_{\Rpn}\lbk\prod_{i\in\Nn\setminus\lbr j\rbr}s_i^{\nu}\rbk \EP\lbk f\lpa\Wx_{\s}\rpa\rbk\prod_{i=1}^{n} \KEmitzs d\s
\\=&\int_{\Rpn}\lbk\prod_{i\in\Nn\setminus\lbr j\rbr}s_i^{\nu}\rbk \partial^{2\nu}_{x_k}\EP\lbk f\lpa\Wx_{\s}\rpa\rbk\prod_{i=1}^{n} \KEmitzs d\s.
\end{split}
\eeq
the mixed partial derivatives cases can be proved similarly.

 Using the boundedness of $f$ and  by observing  the facts
 \beqs
\lim_{y_{k}\to\pm\infty}\pd^{j-1}_{y_{k}}f(y)\pd^{2\nu-j}_{y_{k}}\KBsnjrjxy=0, \ j=1,2, \cdots,  (2\nu-2)
\eeqs
(since all the derivatives up to order $2\nu-3$ are bounded and Lipschitz in $y_k$)  and using integration by parts $(2\nu-2)$-times
we get  (using symmetry of $K^{BS(n,d)}_{s;x,y}$ in $x$ and $y$)
\beq\lbl{outexpec-Lambda}
\bsp
\pd^{2\nu}_{x_{k}}\EP\lbk f(\Wx_{\s})\rbk&=\int_{\Rd} f(y) \pd^{2\nu}_{x_{k}} K^{BS(n,d)}_{s;x,y} dy\\
&=\int_{\Rd} f(y) \pd^{2\nu}_{y_{k}}K^{BS(n,d)}_{s;x,y} dy\\
&=\int_{\Rd}\pd^{2\nu-2}_{y_{k}}f(y) \pd^{2}_{y_{k}}K^{BS(n,d)}_{s;x,y} dy\\
&=\EP\lbk \bigg( \frac{(x_k-W_k^{x_k}(\s))^2-\prod_{i=1}^n s_i}{\prod_{i=1}^n s_i^2}\bigg)(h\lpa\Wx_{\s}\rpa-h(x))\rbk
\end{split}
\eeq
Where $ \pd^{2\nu-2}_{y_{k}}f(y)=h(y)$ is H\"older continuous  with exponent $0<\alpha\leq 1$. Hence  we can finish the proof of part (b) as the last lines of the proof of part (a).
Now
\begin{equation}
\begin{split}
&\int_{\Rpn}\lbk\prod_{i\in\Nn\setminus\lbr j\rbr}s_i^\nu\rbk \bigg|\partial^{2\nu}_{x_k}\EP\lbk f\lpa\Wx_{\s}\rpa\rbk\bigg|\prod_{i=1}^{n}
\KEmitzs d\s\\
&\leq \int_{\Rpn}\lbk\prod_{i\in\Nn\setminus\lbr j\rbr}s_i^\nu\rbk \frac{C}{\prod_{i=1}^n s_i^{1-\alpha/2}}\prod_{i=1}^{n}
\KEmitzs d\s\\
&=C\int_{\Rpn}\lbk\prod_{i\in\Nn\setminus\lbr j\rbr}s_i^{\nu-1+\alpha/2}\rbk s_j^{-(1-\alpha/2)} \prod_{i=1}^{n}
\KEmitzs d\s\\
&=C\bigg(\int_{0}^\infty s_j^{-1+\alpha/2}\KEmitzs ds_{j}\bigg) \prod_{i\in\Nn\setminus\lbr j\rbr}\int_{0}^\infty s_i^{\nu-1+\alpha/2}\KEmitzs d s_{i}\\
&=C \EP[\Lambda(t_j)^{-(1-\alpha/2)}]\prod_{i\in\Nn\setminus\lbr j\rbr}\EP[\Lambda(t_i)^{(\nu-1+\alpha/2)}]\\
&=C t_j^{\lpa\frac\alpha2-1\rpa/\nu}E\lpa\frac1\nu, \frac\alpha2-1\rpa\prod_{i\in\Nn\setminus\lbr j\rbr} t_i^{(\nu-1+\frac\alpha2)/\nu}E\lpa\frac1\nu, \lpa\nu-1+\frac\alpha2\rpa\rpa<\infty.
\end{split}
\end{equation}
hence we obtain \eqref{k-derivative-nu-order} by a standard classical argument.
\epf
\thmref{ISLTBSPDE} can now be proved.
 \bpfs{Proof  of \thmref{ISLTBSPDE}}
Throughout the proof $j\in\Nn$ is fixed but arbitrary.  Let $0<\beta<1$.  Let $v_{\snj}(\rho_{j},x)$ be given by \eqref{vsnjrj}.
Following the same steps of our proof of \thmref{AN-BTBSPDE} above gives
\beq\label{reveallevy}
\bsp
u(\t,x)&=\EP\lbk f\lpa\WxEnt\rpa\rbk
=\int_{\Rpn}\EP\lbk f(\Wx_{\s})\rbk\prod_{i=1}^{n} \KEmitzs d\s
\\&=f(x)+\frac{1}{2}\lbr\int_{\Rpn}\lbk\prod_{i\in\Nn\setminus\lbr j\rbr}s_i\rbk
\int_{0}^{s_{j}}\Delta_{x}v_{\snj}(\rho_{j},x)d\rho_{j}\prod_{i=1}^{n} \KEmitzs d\s\rbr
\end{split}
\eeq
 we then take the $j$-th time Caputo  $\beta$-derivative $\partial^{\beta}_{t_{j}}$, using the dominated convergence theorem to bring $\partial^{\beta}_{t_{j}}$ under the integral; we then use the fact that $\partial_{t_{j}}^{\beta}\KEmjtzs=-\partial_{s_j}\KEmjtzs-\frac{t_{j}^{-\beta}}{\Gamma(1-\beta)}\delta_0(s_j)$\footnote{$\delta_{0}(s_{j})=\lim_{tj\searrow0}\KEmjtzs$ in the sense of tempered distributions ($\delta_{0}(s_{j})ds_{j}$ is the Dirac measure at $0$).} together with integration by parts, observing that  $\int_{0}^{\infty}\int_{0}^{s_{j}}\Delta_{x}v_{\snj}(\rho_{j},x)d\rho_{j}\delta_{0}(s_{j})ds_{j}=0$, along with the fundamental theorem of calculus and \lemref{lapintISLTBS} (a), to get
\beq\lbl{betader}
\bsp
\partial_{t_{j}}^{\beta}u(\t,x)&=\frac{1}{2}\lbr \int_{\Rpn}\lbk\prod_{i\in\Nn\setminus\lbr j\rbr}s_i\rbk \int_{0}^{s_{j}}\Delta_{x}v_{\snj}(\rho_{j},x)d\rho_{j}\pa^{\beta}_{t_{j}} \prod_{i=1}^{n}\KEmjtzs d\s\rbr
\\&=\frac{\Delta_{x}}{2}\lbr \int_{\Rpn}\lbk\prod_{i\in\Nn\setminus\lbr j\rbr}s_i\rbk \EP\lbk f\lpa\Wx_{\s}\rpa\rbk \prod_{i=1}^{n}\KEmjtzs d\s \rbr
\\&=\frac{1}{{2}}\Delta_{x}\sVj(\t,x)
\end{split}
\eeq
The temporal-boundary conditions in \eqref{isstbsfrcsystem}  follow easily from the definitions of $\lbr u,\sVj\rbr$, and the proof of our ISLTBS connection to the fractional PDEs system in \eqref{isstbsfrcsystem} is complete since $j$ is arbitrary.

Now, fix an arbitrary $\beta=1/\nu$ such that $\nu\in\lbr2,3,4,\ldots\rbr$.  Define the simplex
$$\S_{j}^{\nu}:=\lbr \lpa\tau_{1},\tau_{2},...,\tau_{\nu}\rpa\in[0,s_{j}]^{\nu};0\leq \tau_{1}\leq \tau_{2}\leq ...\leq \tau_{\nu}\leq s_{j}\rbr.$$
We use the notation
$$\int_{\S_{j}^{\nu}}h(\tau_{1})\bigotimes_{k=1}^{\nu}d\tau_{k}=\int_{0}^{s_{j}}\cdots\int_{0}^{\tau_{2}}h(\tau_{1})d\tau_{1}\cdots d\tau_{\nu}$$
Let $v_{\snj}\lpa\tau_{1},x\rpa:=\EP\lbk f\lpa\Wx_{\snj}\lpa\tau_{1}\rpa\rpa\rbk$ (the $j$-th temporal parameter is $\tau_{1}$).  Applying It\^o's formula $\nu$ times to $f\lpa W_{\s}^{x}\rpa$ in \eqref{BTBSexp0}, in the ${j}$-th temporal parameter, using the independence of the coordinates of the Brownian-sheet $\Wx$, the fact that the quadratic variation $\lqv\Wx_{\snj}(\cdot)\rqv_{r_{j}}=r_{j}\prod_{i\in\Nn\setminus\lbr j\rbr}s_i$, and taking expectations,  along with repeated use of \lemref{lapintBTBS} (a) and (b), we get
\beq\label{reveallevy2}
\bsp
u(\t,x)&=f(x)+\sum_{k=1}^{\nu-1}\frac{\Delta_{x}^{k}f(x)}{2^{k}}\lbr\int_{\Rpn}\lbk\prod_{i\in\Nn\setminus\lbr j\rbr}s_i\rbk^k \frac{s_{j}^{k}}{k!}\prod_{i=1}^{n} \KEmitzs d\s\rbr
\\&+\frac{1}{2^{\nu}}\lbr\int_{\Rpn}\lbk\prod_{i\in\Nn\setminus\lbr j\rbr}s_i^{\nu}\rbk \int_{\S_{j}^{\nu}}\Delta^{\nu}_{x}v_{\snj}\lpa\tau_{1},x\rpa\bigotimes_{k=1}^{\nu}d\tau_{k}\prod_{i=1}^{n} \KEmitzs d\s\rbr
\end{split}
\eeq
Taking the $j$-th time derivative $\partial_{t_{j}}$, with respect to $t_{j}$, in \eqref{reveallevy2} and the dominated convergence theorem to bring $\partial_{t_{j}}$ under the integral over $\Rpn$; and then using \lemref{lambdapde} together with $\nu$ integration by parts (in which the boundary terms vanish as in the BTBS case---using similar argument as the one before \eqref{fourthBTBS} in the BTBS case above, $\Delta^{\nu}_{x}v_{\snj}\lpa\tau_{1},x\rpa=\Delta_{x}\EP\lbk\Delta^{\nu-1} f\lpa\Wx_{\snj}\lpa\tau_{1}\rpa\rpa\rbk$ is continuous by the conditions on $f$ and its derivatives of orders up to $2\nu-2$  and the conclusion follows) and $\nu$ applications of the fundamental theorem of calculus, along with \lemref{lapintISLTBS} (b), to get
\beq\lbl{highorderISLTBS}
\bsp
\partial_{t_{j}}u(\t,x)&=\sum_{k=1}^{\nu-1}\frac{\Delta_{x}^{k}f(x)}{2^{k}}\partial_{t_{j}}\lbr\int_{\Rpn}\lpa\prod_{i\in\Nn\setminus\lbr j\rbr}s_i\rpa^k \frac{s_{j}^{k}}{k!}\prod_{i=1}^{n} \KEmitzs d\s\rbr
\\&+\frac{\Delta^{\nu}_{x}}{2^{\nu}}\lbr\int_{\Rpn}\lpa\prod_{i\in\Nn\setminus\lbr j\rbr}s_i\rpa^{\nu} \EP\lbk f\lpa\Wx_{\s}\rpa\rbk\prod_{i=1}^{n} \KEmitzs d\s\rbr
\end{split}
\eeq
The temporal-boundary conditions in \eqref{iltbsystem-2m}  follow easily from the definitions of $\lbr u,\sUj_{\nu}\rbr$, and we are done since $j$ is arbitrary.
The terms $\EP[\Lambda(t)^{\gamma}]$ for $-1<\gamma$ as follows:
  First $\Lambda(t)\stackrel{(d)}{=} (L(1)/t)^{-\beta}$ by Corollary 3.1 Meerschaert and Scheffler (2004) hence
  $$
  \EP[\Lambda(t)^{\gamma}]=\EP[(L(1)/t)^{-\gamma\beta}] =t^{\beta\gamma}\EP [L(1)^{-\gamma\beta}]=:t^{\beta\gamma}U(\beta, \gamma)<\infty.
  $$
  Hence
  \begin{equation}
  \begin{split}
  M_\kappa^{(j)}(\t)&=\lbk \prod_{i\in\Nn\setminus\lbr j\rbr}\EP\lpa\Lambda^{(i)}(t_{i})\rpa^\kappa\rbk\frac{\EP\lpa\Lambda^{(j)}(t_{j})\rpa^{\kappa}}{\kappa!}\\
&=\int_{\Rpn}\lbk\prod_{i\in\Nn\setminus\lbr j\rbr}s_i\rbk^{\kappa} \frac{s_{j}^{\kappa}}{\kappa!}\prod_{i=1}^{n} \KEmitzs d\s\\
    &=\frac{U(1/\nu,\kappa)^n}{\kappa !}\prod_{i\in\Nn}t_i^{ \kappa/\nu},\\
  N_\nu^{(j)}(\t)&=\prod_{i\in\Nn\setminus\lbr j\rbr}\EP\lpa\Lambda^{(i)}(t_{i})\rpa^{\nu}\\
  &=\prod_{i\in\Nn\setminus\lbr j\rbr} t^{}U(1/\nu, \nu)=\bigg(\EP[L(1)^{-1}]\bigg)^{n-1}\prod_{i\in\Nn\setminus\lbr j\rbr} t_i,
  \end{split}
  \end{equation}
  completing the proof.
 \epfs
 \subsection{Proof of the conditional equivalence result: the ISLTBS case}
 The following will be useful in the proof of \thmref{equivisltbs}.
\blm\label{iterated-frac-der-lm}
Let $0<\beta_1, \beta_2<1$ and $\beta_1+\beta_2\leq 1$ then we have
$$
\partial_t^{\beta_1}\bigg(\partial_t^{\beta_2} f(t)\bigg)=\partial_t^{\beta_1+\beta_2}f(t)-\frac{t^{\beta_1}}{\Gamma(1-\beta_1)}\bigg(\partial_t^{\beta_2} f(t)\bigg|_{t=0}\bigg),
$$
if all the derivatives exist and are continuous.  If  $\beta=1/\nu$,  $\nu\in\lbr2,3,4,\ldots\rbr$, then
\beq
\begin{split}
\pa_{t}^{\nu\otimes\beta}f(t)&=\partial_tf(t) -\sum_{\kappa=1}^{\nu-1}\frac{t^{-\frac{\kappa}{\nu}}}{\Gamma(1-\frac\kappa\nu)}\bigg[\pa_{t}^{(\nu-\kappa)\otimes\beta}f(t)\bigg|_{t=0}\bigg],
\end{split}
\eeq
if all the derivatives exist and are continuous.
\elm
\bpf
This follows by a simple application of the Laplace transform of the Caputo fractional derivative.
\epf
\brm
We can use Lemma \ref{iterated-frac-der-lm} and the method in the appendix to show that the system $\lbr u,\sUj_{\nu}\rbr$ in \thmref{ISLTBSPDE} satisfy \eqref{iltbsystem-2m}.  First note that
\beq\begin{split}
\bigg[\pa_{t_j}^{k\otimes1/\nu}u\bigg|_{t_j=0}\bigg]&=\frac{\lap^k f(x)}{2^k}\lbk \prod_{i\in\Nn\setminus\lbr j\rbr}\EP\lpa\Lambda^{(i)}(t_{i})\rpa^\kappa\rbk\\
&=\frac{\lap^k f(x)}{2^k}\lbk E(1/\nu,k)^{n-1}\prod_{i\in\Nn\setminus\lbr j\rbr}t_i^{k/\nu}\rbk.
\end{split}
\eeq
Hence we can show that
\beq\label{iltbsystem-2m-form-2}
\begin{split}
&\partial_{t_j}u =\pa_{t_j}^{\nu\otimes\frac1\nu} u+\sum_{k=1}^{\nu-1}\frac{t_j^{-(\nu-k)/\nu}}{\Gamma(1-(\nu-k)/\nu)}
\bigg[\pa_{t_j}^{k\otimes\frac1\nu}u\bigg|_{t_j=0}\bigg]\\
&=\frac{1}{2^{\nu}}\Delta^{\nu}_{x}\sU_{\nu}^{(j)}+\sum_{k=1}^{\nu-1}\frac{t_j^{-(\nu-k)/\nu}}{\Gamma\lpa1-\frac{\nu-k}{\nu}\rpa}
\frac{\lap^k f(x)}{2^k}\lbk E(1/\nu,k)^{n-1}\prod_{i\in\Nn\setminus\lbr j\rbr}t_i^{k/\nu}\rbk
 \end{split}
\eeq
Now comparing the equation \eqref{iltbsystem-2m} and \eqref{iltbsystem-2m-form-2} we see that
$$
E(1/\nu, k)=\EP(\Lambda(1)^k)=\frac{\nu(k-1)!}{\Gamma(k/\nu)}.
$$
\erm
  \bpfs{Proof of Theorem \ref{equivisltbs}}
Fix an arbitrary $j\in\Nn$ and suppose that two pairs $\lbr u,\sVj\rbr$ and $\lbr u,\sUj_{\nu}\rbr$ satisfy \eqref{iltbsystem-2m} and \eqref{isstbsfrcsystem}  for $\beta=1/\nu$.
Applying  $\p_{t_j}^{(\nu-1)\otimes1/\nu}$ to both sides of \eqref{isstbsfrcsystem} and using condition \eqref{u-coefficients} and \lemref{iterated-frac-der-lm} give
\beq \label{iterated-fractional-der}
\begin{split}
&\p_{t_j}^{(\nu-1)\otimes1/\nu} \lbk\frac12\lap_x\mathscr{V}^{(j)}\rbk
\\&=\p_{t_j}^{\nu\otimes1/\nu}u=\partial_{t_j}u-\sum_{k=1}^{\nu-1}\frac{t_j^{-(\nu-k)/\nu}}{\Gamma(1-(\nu-k)/\nu)}
\bigg[\p_{t_j}^{k\otimes1/\nu}u\bigg|_{t_j=0}\bigg]\\
&=\frac{1}{2^{\nu}}\Delta^{\nu}_{x}\sU_{\nu}^{(j)}+\sum_{k=1}^{\nu-1}t_j^{-(\nu-k)/\nu}\frac{ E(1/\nu,k)^n\lap_x^{k}f(x)}{\nu2^{k}(k-1) !}
\\&\quad -\sum_{k=1}^{\nu-1}\frac{t_j^{-(\nu-k)/\nu}}{\Gamma(1-(\nu-k)/\nu)}
\bigg[\p_{t_j}^{k\otimes1/\nu}u\bigg|_{t_j=0}\bigg]
\prod_{i\in\Nn\setminus\lbr j\rbr}t_i^{ \frac k\nu}
\\&=\frac{1}{2^{\nu}}\Delta^{\nu}_{x}\sU_{\nu}^{(j)}
 \end{split}
\eeq
establishing \eqref{equivalence-cond-iltbs} (a) and thus  \eqref{equivalence-cond-iltbs} since \eqref{equivalence-cond-iltbs} (b), (c), and (d) are satisfied by assumption.

For the rest of the proof, suppose that \eqref{u-coefficients} and \eqref{equivalence-cond-iltbs} hold.  We then assume first that $\lbr u,\sVj\rbr$ is a solution to \eqref{isstbsfrcsystem} for $\beta=1/\nu$.
  Applying  $\p_{t_j}^{(\nu-1)\otimes1/\nu}$ to both sides of \eqref{isstbsfrcsystem}  as in
\eqref{iterated-fractional-der}  and using the conditions in  \eqref{u-coefficients} and \eqref{equivalence-cond-iltbs} we get the fact that
 $\lbr u,\sUj\rbr$ is a solution to \eqref{iltbsystem-2m}.

 Now, let $A^{j}(\t,x):=\frac12\lap_x\mathscr{V}^{(j)}(\t,x)$.  If we assume that $\lbr u,\sUj_{\nu}\rbr$ is a solution to \eqref{iltbsystem-2m}, then we can use  \eqref{u-coefficients} and \eqref{equivalence-cond-iltbs} to rewrite \eqref{iltbsystem-2m} as
\beq \label{u-v-pde-m}
\begin{split}
 \p_{t_j}u&=\frac{1}{2^{\nu}}\Delta^{\nu}_{x}\sU_{\nu}^{(j)}+\sum_{k=1}^{\nu-1}t_j^{-(\nu-k)/\nu}\frac{ E(1/\nu,k)^n\lap_x^{k}f(x)}{\nu2^{k}(k-1) !}
\prod_{i\in\Nn\setminus\lbr j\rbr}t_i^{ \frac k\nu}\\
&=\p_{t_j}^{(\nu-1)\otimes1/\nu} A^j+\sum_{k=1}^{\nu-1}\frac{t_j^{-(\nu-k)/\nu}}{\Gamma(1-(\nu-k)/\nu)}
\bigg[\p_{t_j}^{k\otimes1/\nu}u\bigg|_{t_j=0}\bigg]
 \end{split}
\eeq
 Using  Lemma \ref{iterated-frac-der-lm} we get
 \beq \label{u-v-pde-m-ff}
\begin{split}
 \p_{t_{j}}^{\nu\otimes1/\nu}u=\p_{t_j}^{(\nu-1)\otimes1/\nu} A^j.
 \end{split}
\eeq
Now, taking the Fourier-Laplace transform ($x\to k$, and $t_j\to s_j$) yields
\beq
 \begin{split}
 s_j^{1/\nu}\bar{u}_{\tnj}(s_j,k)-s_j^{1/\nu-1}\hat{f}(k)&=   \bar{A}_{\tnj}^{j}(s_j,k)
 \end{split}
\eeq
Taking inverse Fourier-Laplace transform implies that  $\lbr u,\sVj\rbr$ is a solution to \eqref{isstbsfrcsystem} for $\beta=1/\nu$.
\epfs
\section{On relaxing the boundedness condition on $f$ and its derivatives}\lbl{unbdd}As with some problems of applied mathematics, the boundedness assumption on $f$ and its first few derivatives can sometimes be too restrictive.  We now briefly discuss how to relax this boundedness condition.  Carefully examining the proofs of the main results reveal what types of unbounded $f$ are allowed.  We need \lemref{lapintBTBS}  to hold for the BTBS PDEs results to hold, and we need \lemref{lapintBTBS} (a) and (b) and \lemref{lapintISLTBS} for the ISLTBS PDEs results to hold.

In addition to the H\"older continuity on $f$ and all of its $2\nu-2$ order derivatives, $\nu=\frac1\beta\in\lbr2,3,\ldots\rbr$, we need some growth and integrability conditions with respect to the BS density given by \eqref{bsdensity}.  Namely, we can replace the boundedness assumptions on $f$ and its first $2\nu-2$ derivatives by
 \beq\lbl{unbddcond}
\bsp
&\mathrm{(i)}\int_{\Rd}\lab \df{\pd^{i}f(y)}{\pd y^{i_{1}}_{k}\pd y^{i_{2}}_{l}}\df{\pd^{j}\KBsnjrhojxy}{\pd y^{j_{1}}_{k}\pd y^{j_{2}}_{l}} \rab dy<\infty;\ i=0,1,\ldots,2\nu-2, j=0, 1, \ldots,2\nu,\\
&\hspace{5.8cm}\mbox{ and }i+j\le 2\nu,\mbox{ with } k,l=1,\ldots,d,\\
&\mathrm{(ii)}\int_{0}^{r_{j}}\lab\pd^{2}_{x_{l}}\int_{\Rd}\sum_{k=1} ^{d}\pd^{2}_{y_{k}}f(y)\KBsnjrjxy dy\rab d\rho_{j}<\infty;\ l=1,\ldots,d,\ r_{j}>0.
\end{split}
\eeq
Of course, in particular, \eqref{unbddcond} implies
\beq\lbl{unbddgrwth}
\lim_{y_{k}\to\pm\infty}\pd^{j-1}_{y_{k}}f(y)\pd^{2\nu-j}_{y_{k}}\KBssxy=0; \ j=1,2, \cdots,  (2\nu-2),\ k=1,\ldots,d.
\eeq
By the exponential decay in the BS density and its spatial derivatives, condition \eqref{unbddcond} holds when $f$ and all of its derivatives of order up to $2\nu-2$ have polynomial growth.

The interested reader can check that \lemref{lapintBTBS} and \lemref{lapintISLTBS} hold if we     preserve all the differentiability and H\"older continuity conditions on $f$, while replacing the boundedness on $f$ and its derivatives by \eqref{unbddcond}.  This means that, without significant changes in our proofs above (except obviously using  condition \eqref{unbddcond} on $f$ and its derivatives in place of the boundedness condition on $f$ and its derivatives of order  up to $2\nu-2$), we have stronger versions of all our main results that we capture in the following theorem.
\bfr
\bthm
Assume that $\beta=1/\nu$$,$ $\nu\in\lbr2,3,4,\ldots\rbr$.  Then,  \thmref{A-BTBSPDE}--\thmref{equivisltbs} all hold with the boundedness conditions on $f$ and its derivatives of order up to $2\nu-2$ replaced by the conditions in \eqref{unbddcond}, leaving all other conditions in force.
\ethm
\efr
We note that the BTBS results (\thmref{A-BTBSPDE}, \thmref{AN-BTBSPDE}, and \thmref{equiv}) fall under the case of $\nu=1/2$, with a minor scaling of the density ($4t$ is replaced by $2t$), as in equation \eqref{Ydens}.

\subsection*{Acknowledgement}  We are sincerely grateful for the anonymous referee's careful and constructive comments which improved the readability of the paper.  In particular, Subsection \ref{fexample} and Section \ref{unbdd} were added as an answer to an excellent point made by the referee.
 \appendix
 \section{Another proof of the fractional PDEs connection in \thmref{ISLTBSPDE}}
Here, we present a different proof of the fractional PDEs connection \thmref{ISLTBSPDE} that is an adaptation of the proof of Theorem 1.1 in \cite{Abtbs}.  Take the time fractional derivative and put it inside the integral by the dominated convergence theorem,
 using part (e) of Lemma \ref{fractional-density-pde} and integration by parts once, and the boundary conditions $\lim_{s_j\to\infty}\KEmjtzs=0$ and $\lim_{s_j\to0}\KEmjtzs=-\frac{t_{j}^{-\beta}}{\Gamma(1-\beta)}$, and using Lemma 1.1 in \cite{Abtbs} to get
\begin{equation}\begin{split}
&\partial_{t_j}^\beta u(t,x)
=\int_{\Rp^n}\EP \lbk  f\lpa\Wpsx\rpa\rbk\ \partial_{t_j}^\beta\KEmjtzs \prod_{i\in\Nn\setminus\lbr j\rbr} \KEmitzs d\s \\
&=\int_{\Rp^n}\EP\lbk f\lpa\Wpsx\rpa\rbk\ \lpa-\partial_{s_j}\KEmjtzs-\frac{t_{j}^{-\beta}}{\Gamma(1-\beta)}\delta_0(s_j) \rpa\prod_{ i\in\Nn\setminus\lbr j \rbr} \KEmitzs d\s\\
&=\int_{\Rp^n}\partial_{s_j} \EP \lbk f \lpa \Wpsx\rpa \rbk\ \prod_{i=1}^n \KEmitzs d\s\\
& +\int_{\Rp^{n-1}} f \lpa x\rpa \lpa-\frac{t_{j}^{-\beta}}{\Gamma(1-\beta)}\rpa \prod_{i\in\Nn\setminus\lbr j\rbr} \KEmitzs d\snj \\
&+\int_{\Rp^{n-1}} \EP\lbk f \lpa\Wpsx\rpa  \rbk \KEmjtzs \bigg|_{s_j=0} \prod_{i\in\Nn\setminus\lbr j\rbr} \KEmitzs d\snj\\
&= \lap_x  \int_{\Rp ^n}\lpa \frac12 \prod_{i\in\Nn\setminus \{ j \} } s_i \rpa \EP\lbk f\lpa\Wpsx\rpa\rbk\ \prod_{i=1}^n \KEmitzs d\s =\frac12 \lap_x \mathscr{V}^{(j)}(\t,x).\\
\end{split}
\end{equation}
To justify taking the Laplacian outside of the integral we use \lemref{lapintISLTBS}.
Next we verify the boundary condition (d), since at $t_j=0$, $ \EP\lbk f\lpa\Wpsx\rpa\rbk=f(x)$ we get
\begin{equation}
\begin{split}
\mathscr{V}^{(j)}(\t,x)&=\int_{\Rp ^n}\lpa  \prod_{i\in\Nn\setminus \{ j \} } s_i \rpa f(x) \prod_{i=1}^n \KEmitzs d\s\\
&=\int_{\Rp ^{n-1}}\lpa  \prod_{i\in\Nn\setminus \{ j \} } s_i \rpa f(x) \prod_{i\in\Nn\setminus \{ j \}} \KEmitzs d\s\\
&=f(x)\prod_{i\in\Nn\setminus \{ j \}}   \EP \lpa \Lambda(t_i)\rpa \\
&=f(x)\prod_{i\in\Nn\setminus \{ j \}}\lpa C(\beta)t_i^\beta \rpa\\
\end{split}
\end{equation}
Where the last equality follows from Corollary 3.1 in \cite{limitCTRW}
\section{Frequent acronyms and notations key}\lbl{B}
\begin{enumerate}\renewcommand{\labelenumi}{\Roman{enumi}.}
\item {\textbf{Acronyms}}\vspace{2mm}
\begin{enumerate}\renewcommand{\labelenumii}{(\roman{enumii})}
\item BTBM: Brownian-time Brownian motion.
\item BTBS: Brownian-time Brownian sheet.
\item BTP: Brownian-time process.
\item ISLTBS: inverse-stable-L\'evy-time Brownian sheet
\end{enumerate}
\vspace{2.5mm}
\item {\textbf{Notations}}\vspace{2mm}
\begin{enumerate}\renewcommand{\labelenumii}{(\roman{enumii})}
 \item $\Nn=\lbr 1,\ldots,n\rbr$.
\item $\KBmitzs$:  The density of a one-dimensional BM starting at $0$ (see Section \ref{inverse} just before \eqref{Ydens}).
\item $\KEmitzs$:  The density of a one-dimensional  $\beta$-inverse L\'evy motion $\Lambda$ astarting at $0$ (see Section \ref{inverse}).
\item $\KBstxy$: The density (or kernel) of an $n$-parameter $d$-dimensional Brownian sheet (see \eqref{bsdensity}).
\item  $\Delta_{x}$ vs.~$\Delta$ (see Remark \ref{laplacians})\vspace{0.5 mm}
\end{enumerate}
\end{enumerate}

\end{document}